\documentclass[12pt]{article}
\usepackage[utf8]{inputenc}
\usepackage{times}
\usepackage{graphicx}
\usepackage{amssymb}
\usepackage{amsthm}
\usepackage{url,hyperref}
\usepackage{algorithm}
\usepackage{algorithmic}
\usepackage{amsfonts}
\usepackage{indentfirst} 
\usepackage{natbib}
\usepackage{mathrsfs}
\usepackage{multirow}
\usepackage{comment}
\usepackage{bm}
\usepackage{booktabs} 
\usepackage{float}
\usepackage{tikz}
\usepackage{amsmath}
\usepackage{threeparttable}
\usepackage{setspace}
\usepackage{enumitem}
\usepackage{subfig} 
\usepackage{pdfpages}
\usepackage[modulo,displaymath,mathlines]{lineno}
\allowdisplaybreaks

\usepackage{xr}
\makeatletter
\newcommand*{\addFileDependency}[1]{% argument=file name and extension
	\typeout{(#1)}
	\@addtofilelist{#1}
	\IfFileExists{#1}{}{\typeout{No file #1.}}
}
\makeatother

\newcommand*{\myexternaldocument}[1]{%
	\externaldocument{#1}%
	\addFileDependency{#1.tex}%
	\addFileDependency{#1.aux}%
}
\myexternaldocument{supplement}

\def\spacingset#1{\renewcommand{\baselinestretch}%
{#1}\small\normalsize} \spacingset{1}

\usepackage{hyperref}
\hypersetup{
    citecolor={blue!70!black},
    colorlinks=true,
    linkcolor={red!50!black},
    filecolor={red!50!black},
}

\setlist[description]{%
  font={\normalfont}, % set the label font
}

\usepackage{titlesec}
\titleformat*{\paragraph}{\itshape}

\def\lin#1{#1}

\def\real{\mathbb{R}}

\def\tdomain{\mathcal{T}}
\def\batchsize{B}
\def\sPhiPhi{H}
\def\sPhy{G}

\def\expect{\mathbb{E}}
\def\prob{\mathrm{Pr}}
\def\Op{O_{P}}
\def\asympleq{\lesssim}
\def\asympgeq{\gtrsim}

\newtheorem{theorem}{Theorem}
\newtheorem{proposition}{Proposition}

\newtheorem{corollary}{Corollary}
\newtheorem{assumption}{Assumption}

\newtheorem{example}{Example}

\newcounter{remark}
\newenvironment{remark}[1][]{\refstepcounter{remark}\par\medskip
   \noindent \textit{Remark~\theremark. #1} }{\medskip}

\makeatletter
  \long\def\myempty{}
  \def\XR@addURL#1{\XR@@dURL#1\myempty{}{}{}{}{}\\}
  \def\XR@@dURL#1#2#3#4#5#6#7\\{%
    {#1}{#2}%
    \ifx\myempty#6\@empty
      {#3}{#4}{\XR@URL}%
    \else
    \fi
  }
\makeatother

\graphicspath{{fig/}}

\def\suppref#1#2{#2}

\date{}

\def\mytitle{Optimal One-pass Nonparametric Estimation Under Memory Constraint}

\newcommand{\blind}{1}

\begin{document}

\nolinenumbers

\def\spacingset#1{\renewcommand{\baselinestretch}%
{#1}\small\normalsize} 
\spacingset{1}

\if1\blind
{
  \title{\bf \mytitle}
  \author{Mingxue Quan\thanks{
    Mingxue Quan (quanmingxue@ruc.edu.cn) was a visiting Ph.D. student in National University of Singapore when conducting the research.  Research is supported by Chinese Government Scholarship (CSC202006360089).}\hspace{.2cm}\\
    \textit{School of Mathematics, Renmin University of China}\\
    ~\\
    and \\
    ~\\
    Zhenhua Lin\thanks{Correspondence: linz@nus.edu.sg. Research is partially supported by NUS Startup Grant A-0004816-01-00.} \\
    \textit{Department of Statistics and Data Science, National University of Singapore}}
  \maketitle
} \fi

\if0\blind
{
  \title{\bf \mytitle}
  \author{}
  \bigskip
  \bigskip
  \bigskip
  %\begin{center}
  %  {\LARGE\bf \mytitle}
  %\end{center}
  \maketitle
  \medskip
} \fi

\bigskip
\begin{abstract}
    For nonparametric regression in the streaming setting, where data constantly flow in and require real-time analysis, a main challenge is that data are cleared from the computer system once processed due to limited computer memory and storage. We tackle the challenge by proposing a novel one-pass estimator based on penalized orthogonal basis expansions and developing a general framework to study the interplay between statistical efficiency and memory consumption of estimators. We show that, the proposed estimator is statistically optimal under memory constraint, and has asymptotically minimal memory footprints among all one-pass estimators of the same estimation quality. Numerical studies demonstrate that the proposed one-pass estimator is nearly as efficient as its non-streaming counterpart that has access to all historical data. 
\end{abstract}

\noindent%
    
\vspace{0.5cm}    
\noindent    
\textit{MSC2020 subject classification}: 62G05, 62G08.

\newpage
\spacingset{1.9} % DON'T change the spacing!

%\linenumbers
%\resetlinenumber

\section{Introduction}\label{sec:introduction}
Nonparametric regression, in which  a response variable is paired with a predictor without a prescribed form, is a fundamental statistical tool for data analysis. It has been extensively studied in the literature, among many others, by \cite{nadaraya1964,watson1964,fan1993local} on local smoothing methods, by \cite{hastie1990,wahba1990,green1994} on spline methods and by \cite{stone1982optimal,Donoho1998} on minimax estimation. %For a more comprehensive  treatment on this subject, we recommend the introductory text \cite{Tsybakov2008}. 

Classic nonparametric regression algorithms assume that data of interest have been fully collected at the time of analysis and there is sufficient computer memory to hold the whole data. This assumption is nowadays challenged by the massive volume of modern data for example collected from a large number of wearable devices. In addition, data that constantly flow in and require real-time analysis, like those from wearable devices, also violate the above assumption and call for new analysis paradigms; see Section \ref{subsec:realdata} for a case study on data of this kind. A paradigm particularly suited for real-time and memory-efficient data processing may assume that data points come sequentially and are cleared from memory once processed. Data in such a paradigm are called streaming data or data streams, and algorithms following this paradigm are called streaming or one-pass algorithms.

Nonparametric regression for streaming data has gained increasing  attention in the past two decades. For example, 
\cite{revesz1977,mokkadem2008revisiting} studied recursive kernel estimators, while \cite{vilar2000recursive,gu2012sequential} considered one-pass local polynomial smoothers. %and showed that their estimators achieve the usual minimax $L^2$ convergence rate. 
\cite{huang2013recursive} studied a recursive kernel estimator for time series, while \cite{yang2021online} proposed a one-pass local linear smoother with application to  functional data analysis.
 \cite{yao2021spline} studied penalized spline smoothing for streaming data, focusing on strategies to dynamically place new knots. Although these endeavors consider estimation in the streaming setting, few consider the interplay between statistical efficiency and memory constraint.

The main challenge of nonparametric regression for streaming data is to produce a real-time efficient estimate of the regression function by using limited computer memory while data come sequentially and are not retained once processed. Take for example, the popular basis expansion approach \citep[][among many others]{Chen2015} of which the basic idea is to approximate the regression function by a linear combination of a few basis functions and then estimate the corresponding coefficients from data. In the streaming setting, as data flow in continuously, it is necessary to \emph{dynamically} and periodically add a new basis function into the approximation to improve approximation quality; otherwise, the approximation error, which quantifies the discrepancy between the regression function and the best linear combination, will dominate the estimation error at some time point. An estimate of the coefficient of the new  basis function then needs to be provided at the time when the basis function is added, but the historical data needed to produce the estimate are not retained and thus not available in the streaming setting.

To tackle the above challenge, in this paper we contribute a novel one-pass nonparametric estimation method that  exploits penalized basis expansions and requires only $O(q)$ memory units, where $q$ is the number of basis functions being used. To appreciate this low level of memory footprints and the difficulty to achieve it, one notes that a $q\times q$ matrix $\Phi^\top \Phi$, where $\Phi$ (defined in Section \ref{subsec:model}) is an $n\times q$ matrix containing information about the predictor, is often needed to compute the coefficients, and maintaining this matrix needs at least $q^2$ memory units; see Section \ref{subsec:estimation}. Our strategy for reducing memory consumption is to approximately reproduce the matrix $\Phi^\top\Phi$ at the query time from  a dynamically maintained estimate of the probability density of the predictor by using only $O(q)$ memory units. The proposed method also includes a sophisticated pre-estimation strategy to dynamically and periodically add new basis functions; see Section \ref{subsec:increase-basis-functions}.

Another key contribution is a general framework for investigating the interplay between statistical estimation quality and computational memory consumption, based on which a tight lower bound on the required memory units is established for any one-pass estimator that possesses certain convergence rate in nonparametric regression. Roughly speaking, if the regression function is $\beta$-times continuously differentiable, then  the minimum number of memory units required by any one-pass estimator converging at the rate $\Op(b_n)$ is at least of the order $b_n^{-1/\beta}$. Conversely, if there are $b_n^{-1/\beta}$  available memory units, then the proposed one-pass estimator converges at the rate $\Op(b_n)$. Consequently, \emph{the proposed  estimator is statistically optimal under memory constraint, and computationally optimal in memory consumption among all one-pass estimators of the same estimation quality}. To the best of our knowledge, this is the first general framework and result on the interplay between the nonparametric convergence rate and  memory consumption. It shows that a  more  accurate one-pass estimator generally consumes more memory.  In particular, if the available memory is capped by a universal constant, then there must be a gap, independent of the sample size $n$, between any one-pass estimator and the unknown regression function.
The lower bound is established by leveraging the theory of communication complexity that is well studied in computer science and information theory but less explored in statistics. The key observation is that a nonparametric estimator with quality guarantee may be used by two parties to solve certain problems via a one-way communication channel; see Section \ref{subsec:ccct} for details.

We remark that, the one-pass local  smoother of \cite{yang2021online}, consuming a constant amount of memory, only estimates the regression function at a  finite set of predetermined locations during the streaming process; see Remarks \ref{rem:diff-opls-proposed} and \ref{rem:opls-space-bound} for a detailed discussion on why a constant amount of memory is not sufficient for consistently estimating the entire function in the streaming setting. In addition, compared to our method, the spline smoothing method of \cite{yao2021spline} requires asymptotically larger memory to achieve the same estimation equality.

\section{One-pass Nonparametric Estimation}\label{sec:methodology}
\subsection{Model and Penalized Basis Expansion}\label{subsec:model}
Let $(T_i,Y_i)$, $i=1,2,\ldots$, be streaming data sampled from the  nonparametric regression model
\begin{equation}\label{eq:npreg}
    Y_i=m(T_i)+\varepsilon_i,
\end{equation}
where $m$ is an unknown regression function capturing the relationship between the response $Y_i\in\real$ and the predictor $T_i\in\tdomain$ with $\tdomain$ being a closed interval of $\real$, and $\varepsilon_i$ represents the measurement error with $\expect\varepsilon_i=0$. The goal is to estimate $m(t)$ \emph{simultaneously for all $t\in\tdomain$} from the streaming data pairs $(T_i,Y_i)$, for which we adopt the penalized basis expansion approach.

Specifically, let $\phi_1,\phi_2,\ldots$ be a predetermined sequence of basis functions such that $m(t)=\sum_{j=1}^\infty a_j\phi_j(t)$ for coefficients $a_1,a_2,\ldots$. Examples of such basis functions include the Fourier basis functions when $m$ is periodic or the wavelet basis functions \citep{cohen1993wavelets}; see Section \ref{subsec:basis} for details. 
In practice, only the first $q$ leading basis functions are used, which result in an approximation $\sum_{j=1}^q a_j\phi_j(t)$ to $m(t)$. In a non-streaming setting, the  coefficient $a=(a_1,\ldots,a_q)^\top$ may be  estimated by $\tilde a$ that  minimizes the penalized empirical mean squared error, i.e.,
\begin{equation}\label{eq:aobj}
\tilde{a}=\mathop{\arg\min}_{\check a\in \real^q} \bigg\{ \frac{1}{n}\sum_{i=1}^n [Y_i-\check a^\top \phi(T_i)]^2+\rho \check a^\top W\check a\bigg\},   
\end{equation}
where $\phi(t)=(\phi_1(t),\ldots,\phi_q(t))^\top$, $W$ is a $q\times q$ positive semi-definite matrix imposing a ridge-type penalty, and $\rho$ is a tuning parameter that controls the degree of penalty. Examples of $W$ include the $q\times q$ identity matrix $I_q$ and the matrix formed by the elements 
$\int_{\tdomain}\phi_j^{\prime\prime}(t)\phi_k^{\prime\prime}(t)dt$ for $j,k=1,\ldots,q$; the former leads to the classic Tikhonov regularization, while the later corresponds to the roughness penalty since $\check a^\top W\check a = \int_{\tdomain} [\{\check a^\top \phi(t)\}^{\prime\prime}]^2dt$ precisely measures roughness of the candidate estimate $\check a^\top \phi(t)$.  
The minimization problem \eqref{eq:aobj} is solved by 
\begin{equation}\label{eq:a-full}
 \tilde{a}=(n^{-1}\Phi^\top \Phi+\rho W)^{-1}(n^{-1}\Phi^\top Y)\qquad\text{and}\qquad \tilde{m}(t)=\tilde{a}^\top \phi(t),
\end{equation}
where $Y=(Y_1,\ldots,Y_n)^\top$, $\Phi$ is the $n\times q$ matrix with elements $\Phi_{ij}=\phi_j(T_i)$, and $\tilde m(t)$ is a non-streaming estimate of $m(t)$. % is then provided by 
% \begin{equation}
% \tilde{m}(t)=\tilde{a}^\top \phi(t).
% \end{equation}
In the non-streaming setting where all data are retained, when a new batch of data is received, the estimate $\tilde a$ and hence also the estimate $\tilde m$ can be recomputed according to \eqref{eq:a-full} by using the historical data and the new data together. Unfortunately, such recomputation no longer applies to streaming data.

\subsection{Space-saving One-pass Algorithm}\label{subsec:estimation}
For streaming data, by inspecting the solution \eqref{eq:a-full}, one finds that it is sufficient to dynamically maintain the $q\times q$ matrix $\tilde\sPhiPhi:=n^{-1}\Phi^\top\Phi$ and the vector $\sPhy:=\Phi^\top Y$ of length $q$; the matrix $W$ is determined by the predetermined basis functions and thus can be computed on the fly. When new data arrive, $\tilde\sPhiPhi$ and $\sPhy$ may be updated, as follows. Without loss of generality, assume data arrive or are processed in batches and each batch contains $\batchsize\geq 1$ pairs of observations. Let $\Phi_{[n]}$ and $Y_{[n]}$, computed from the $n$th batch of the data, be respectively the counterparts of $\Phi$ and $Y$. Let $\tilde\sPhiPhi_{1:n}$ and $\sPhy_{1:n}$ be respectively the counterparts of $\tilde\sPhiPhi$ and $\sPhy$  after the first $n$ batches are processed. They are updated upon receiving the $n$th batch of data according to 
\begin{equation}\label{eq:update-H-G}
\tilde\sPhiPhi_{1:n}=\frac{1}{nB}\{B(n-1)\tilde\sPhiPhi_{1:(n-1)}+\Phi_{[n]}^\top \Phi_{[n]}\}\qquad\text{and}\qquad
\sPhy_{1:n}=\sPhy_{1:(n-1)}+\Phi_{[n]}^\top Y_{[n]}
\end{equation}
for $n > 1$, with $
\tilde\sPhiPhi_{1:1}= \frac{1}{B}\Phi_{[1]}^\top \Phi_{[1]}$ and 
$\sPhy_{1:1}=\Phi_{[1]}^\top Y_{[1]}$.
This  na\"ive algorithm requires $\Omega(q^2)$ memory units to store the summary statistics $\tilde \sPhiPhi_{1:n}$ and $G_{1:n}$, where $\Omega(q^2)$ denotes a quantity that is asymptotically at least of the order $q^2$. 

To further reduce memory footprints, we observe that  $\tilde\sPhiPhi_{jl}=n^{-1}\sum_{i=1}^n \Phi_{j}(T_i)\Phi_l(T_i)$, the $(j,l)$-element of the matrix $\tilde\sPhiPhi$, may be approximated by $H_{jl}=\expect \{\phi_{j}(T)\phi_l(T)\}$, since $\tilde\sPhiPhi_{jl}$ is the average of random observations whose mean is $\sPhiPhi_{jl}$. 
When the density $f$ of $T$ is known,  $\expect \{\phi_{j}(T)\phi_l(T)\}$ can be computed by $\int_{\tdomain} \phi_j(t)\phi_l(t)f(t)dt$.
Otherwise, we may dynamically maintain an estimator $\hat f$ for $f$ in the memory. For this purpose, we adopt the orthogonal series density estimation (OSDE) method \citep{hall1986rate}, since it requires a low level of memory footprints while has a high level of accuracy; see Proposition \ref{prop:density-rate-Holder-class}. 
Specifically, for predetermined orthogonal basis functions $\psi_1,\ldots,\psi_p$ with $p\geq 1$, the non-streaming OSDE estimate $\tilde f(t)$ is given by 
%\begin{equation}\label{eq:ftilde}
$\tilde{f}(t)=\sum_{j=1}^p\tilde\theta_j\psi_j(t)$ with $\tilde\theta_j=\frac{1}{n}\sum_{i=1}^n \psi_j(T_i).$
%\end{equation}
For streaming data, let $\hat\theta_{1:n,j}$ be the counterpart of $\tilde\theta_j$ after the first $n$ batches of data are processed, with the convention that $\hat\theta_{1:0,j}=0$. We compute $\hat\theta_{1:n,j}$ and obtain the one-pass estimate $\hat f(t)$ of $f(t)$ upon receiving the $n$th batch by
%\begin{equation}\label{eq:thetaupdate}
%\theta_{1:k,j}=\frac{1}{k}\big\{(k-1)\theta_{1:(k-1),j}+\psi_j(T_i)\big\}.
%\end{equation}
%\quan{Consider $\batchsize$, the estimation of $\theta$ is:
\begin{equation}\label{eq:fhat}
\hat\theta_{1:n,j}=\frac{1}{n\batchsize}\big\{\batchsize(n-1)\hat\theta_{1:(n-1),j}+\sum_{i=n\batchsize-\batchsize+1}^{n\batchsize}\psi_j(T_i)\big\}\qquad\text{and}\qquad \hat f(t)=\sum_{j=1}^p\hat\theta_{1:n,j}\psi_j(t).
\end{equation}
The element $\tilde\sPhiPhi_{jl}$ of the matrix $\tilde\sPhiPhi$ is then approximated by
%\begin{equation}\label{eq:Hhat}
   $\hat \sPhiPhi_{jl}=\int_{\tdomain}\phi_j(t)\phi_l(t)\hat f(t)dt$. 
%\end{equation}
Consequently, the matrix $\tilde\sPhiPhi$ no longer needs to be maintained in the memory.

\begin{remark}
Although the estimated function $\hat f$ may not be integrated to one or always non-negative, Proposition \ref{prop:density-rate-Holder-class} in Section \ref{sec:convergence rate} shows that $\hat f(t)$ rapidly converges to $f(t)$ uniformly, so that $\hat\sPhiPhi_{jl}$ also converges to $\tilde \sPhiPhi_{jl}$. In the application where a genuine density is essential, one may use $\max\{0,\hat f(t)\}/\int_{\tdomain}\max\{0,\hat f(s)\}ds$ as an alternative density estimator which has the same uniform convergence rate of $\hat f(t)$.
\end{remark}
%\vspace{5mm}

Finally, to produce a one-pass estimate $\hat m$ by using the above space-saving idea, let $\hat{\sPhiPhi}_q$ be the $q\times q$ matrix whose $(j,l)$-element is  $\int_{\tdomain}\phi_j(t)\phi_l(t)f(t)dt$  when $f$ is known or $\hat{\sPhiPhi}_{jl}$ when $f$ is estimated from the data. 
The \emph{space-saving one-pass estimator} $\hat m$ at the current time $n$,  graphically illustrated in Figure \suppref{fig:graph-alg}{S1} of the supplement, is then given by \begin{equation}\label{eq:mhat-space-saving}
\hat m(t)=\sum_{j=1}^q \hat a_{j}\phi_j(t)\qquad\text{with}\qquad \hat a=(\hat{\sPhiPhi}_{q}+\rho W)^{-1}(N_q^{-1}G_q),
\end{equation}
where $G_q$ is the vector formed by the first $q$ elements of $G_{1:n}$, and $N_q=\mathrm{diag}(n_1,\ldots,n_q)$ with $n_j=n-\tau_j+1$ for $j=1,\ldots,q$ and $\tau_j$ (defined in Section \ref{subsec:increase-basis-functions}) being the start time of pre-estimating the coefficient $a_j$; the matrix $N_q^{-1}$ is introduced to account for the fact that the elements of the vector $G_{1:n}$ may be computed from different amounts of observations. Selection of the tuning parameter $\rho$ is detailed in Section \ref{subsec:tuning-parameter}. It is worth  noting that there is no need to keep the matrix $\hat\sPhiPhi_q$ in the memory since it can be computed on the fly by using the predetermined basis functions $\phi_1,\ldots,\phi_q$ and the density estimator $\hat f$.

The above algorithm, referred to as the {space-saving one-pass algorithm}, requires only $O(q+p)=O(q)$ memory units to maintain the summary statistics $G_{1:n}$ and $\hat\theta_{1:n}=(\hat\theta_{1:n,1},\ldots,\hat\theta_{1:n,p})^\top$, provided $p=O(q)$. Compared with the na\"ive method that requires $\Omega(q^2)$ memory units, the space-saving algorithm substantially improves memory efficiency. In addition, {the choice of  $q$ can be made to respect potential memory constraint}. For example, when there are only $\mathfrak s$ available memory units, we choose $q=O(\mathfrak s)$. In Section \ref{sec:space-lower-bound} we shall show that the space-saving one-pass estimator is statistically optimal under memory constraint.

\begin{remark}\label{rem:diff-opls-proposed}
\lin{The summary statistics $G_{1:n}$ and $\hat\theta_{1:n}$ are \emph{shared} by {all} $t\in\tdomain$, i.e., when an estimate of $m(t)$ for any $t$ is requested, they are combined with the basis functions to produced the estimate. This feature enables the proposed algorithm to provide an estimate of $m(t)$ simultaneously for all $t\in\tdomain$ in the streaming setting. In contrast, the one-pass local smoothing (OPLS) algorithm of \cite{yang2021online} needs to maintain an \emph{individual} set of local statistics \citep[][Section 3]{yang2021online} pertaining to each $t$ in order to estimate $m(t)$. It is clearly infeasible to maintain local statistics for all $t\in\tdomain$ in finite memory. Consequently, the OPLS algorithm needs to predetermine a finite subset $\tdomain_\ast\subset\tdomain$ before it starts to process data streams, and produces estimates of $m(t)$ only for those  $t\in\tdomain_\ast$ but not for other $t\in\tdomain$ during the streaming process.}
\end{remark}

\subsection{Pre-estimation Strategy}\label{subsec:increase-basis-functions}
As mentioned in the introduction, as data continuously flow in, it is necessary to dynamically increase the number of basis functions in estimation in order to further reduce approximation error. However, at the time when a new basis function is added, there are few data to estimate its coefficient since data are not retained in the streaming setting.
To overcome this hurdle, we propose to pre-estimate the coefficient, as follows.

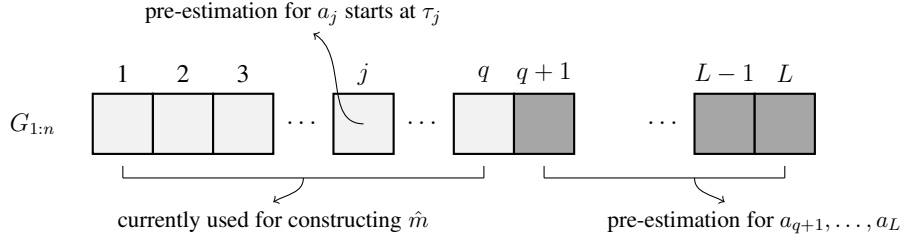
\begin{figure}
    \centering
\begin{tikzpicture}[darkstyle/.style={draw,fill=gray!70,thick},lightstyle/.style={draw,fill=gray!10,thick},scale=0.8,every node/.style={scale=0.8}]
	\def\xa{0}
	\def\k{3}
	\def\a{1}
	
	\node at (0,-\a/2) {$G_{1:n}$};
	
	%\draw (0,0) -- (\k*\a,0);
	%\draw (0,-\a) -- (\k*\a,-\a);
	\foreach \x in {1,...,\k}
	{
		%\draw (\x*\a,0) -- (\x*\a,-\a);
		\filldraw[lightstyle] (\x*\a,0) rectangle ++(\a,-\a);
		\node at (\x*\a+\a/2,\a/3) {\x};
	}
	
	\pgfmathsetmacro{\k}{\k+1}
	\node at (\k*\a+\a/2,-\a/2) {$\cdots$};
	
	\pgfmathsetmacro{\k}{\k+1}
	\pgfmathsetmacro{\j}{\k}
	\filldraw[lightstyle] (\k*\a,0) rectangle ++(\a,-\a);
	\node at (\k*\a+\a/2,\a/3) {$j$};
	
	\pgfmathsetmacro{\k}{\k+1}
	\node at (\k*\a+\a/2,-\a/2) {$\cdots$};
	
	\def\q{2}
	\pgfmathsetmacro{\k}{\k+1}
	\pgfmathsetmacro{\q}{\k+1}
	\filldraw[lightstyle] (\k*\a,0) rectangle ++(\a,-\a);
	
	\pgfmathsetmacro{\k}{\k+1}
	\filldraw[darkstyle] (\k*\a,0) rectangle ++(\a,-\a);
	
	\node at (\q*\a-\a/2,\a/3) {$q$};
	\node at (\q*\a+\a/2,\a/3) {$q+1$};
	
	\pgfmathsetmacro{\k}{\k+2}
	\node at (\k*\a+\a/2,-\a/2) {$\cdots$};
	
	\def\L{2}
	\pgfmathsetmacro{\k}{\k+1}
	\pgfmathsetmacro{\L}{\k+1}
	\foreach \x in {\k,...,\L}
	{
		%\draw (\x*\a,0) -- (\x*\a,-\a);
		\filldraw[darkstyle] (\x*\a,0) rectangle ++(\a,-\a);
	}
	
	\node at (\L*\a-\a/2,\a/3) {$L-1$};
	\node at (\L*\a+\a/2,\a/3) {$L$};
	
	\def\b{\a/5}
	
	\draw (3*\a/2,-\a-\b) -- (3*\a/2,-\a-2*\b);
	\draw (3*\a/2,-\a-2*\b) -- (\q*\a-\a/2,-\a-2*\b);
	\draw (\q*\a-\a/2,-\a-\b) -- (\q*\a-\a/2,-\a-2*\b);

	\draw (\q*\a+\a/2,-\a-\b) -- (\q*\a+\a/2,-\a-2*\b);
	\draw (\q*\a+\a/2,-\a-2*\b) -- (\L*\a+\a/2,-\a-2*\b);
	\draw (\L*\a+\a/2,-\a-\b) -- (\L*\a+\a/2,-\a-2*\b);

	\draw[->] (\j*\a+\a/2,-\a/2) to[out=-180,in=-30] (\j*\a-\a/3,\a);
	\node[above] at (\j*\a-2*\a/3,\a) {\small pre-estimation for $a_j$ starts at $\tau_j$};
	
	\draw[->] (\j*\a-\a/2,-\a-2*\b) to[out=-90,in=30] (\j*\a-\a,-\a-4*\b);
	\node[below] at (\j*\a-\a,-\a-4*\b) {\small currently used for constructing $\hat m$};
	
	\draw[->] (\L*\a-\a,-\a-2*\b) to[out=-90,in=150] (\L*\a-\a/2,-\a-4*\b);
	\node[below] at (\L*\a,-\a-4*\b) {\small pre-estimation for $a_{q+1},\ldots,a_L$};
	
\end{tikzpicture}
\caption{Graphical illustration of the pre-estimation procedure, where each slot represents an element of the summary statistic $G_{1:n}$ and $L$ is the length of $G_{1:n}$ at the timestamp $n$.\label{fig:graph-pre-est}}
\end{figure}

Let $S(j)$ be the time point when the $j$th basis function is activated, i.e., the first time when it is used to produce the estimate $\hat m$, where $S$ is a monotonically increasing function. For example, if $q$ grows with the sample size $n$ at a polynomial rate $n^h$  for some $h\in(0,1)$, then we may choose $S(j)=S_h(j)=\lfloor (C_qj)^{1/h}\rfloor$ for some constant $C_q>0$,
where the notation $\lfloor x\rfloor$ denotes the largest integer that does not exceed $x$; further discussions on the choice of $h$ are provided  in Section \ref{subsec:tuning-parameter}. 
Take a constant $c_\circ\in(0,1)$ such that there exists another constant $c_\star \geq 2$ satisfying $c_{\circ}S(c_\star j)-1 > S(j+1)$ for all $j\geq 1$; such $c_\circ$ exists for the aforementioned function $S_h$.  Now at the time $\tau_j=\lfloor c_\circ S(j)\rfloor$, we start to pre-estimate the coefficient of the $j$th basis function, by appending a zero entry to the vector $G_{1:n}$ and updating it according to \eqref{eq:update-H-G}. Figure \ref{fig:graph-pre-est} provides a graphical illustration of pre-estimation.

According to the above pre-estimation strategy, at the time when the $j$th basis is activated, its corresponding coefficient estimate $\hat a_j$ has been calculated based on at least $n-\tau_j+1 \geq n-c_\circ S(j) + 1\geq (1-c_{\circ})n$ observations since $S(j)\leq S(q)\leq n \leq S(q+1)$ for all $n\geq 1$ and $j=1,\ldots,q$, where we note $q=q_n$ grows with $n$. In addition, the  inequalities $\tau_{c_\star q} \geq c_{\circ}S(c_\star q)-1> S(q+1)\geq n$ suggest that the start time of pre-estimation for the coefficient of the $(c_\star q)$th basis function is later than the current timestamp $n$, and consequently the length of the vector $G_{1:n}$ is less than $c_\star q$. This shows that the amount of memory units required by the pre-estimation strategy is still $O(q)$ for all $n$. The same strategy is adopted to maintain the summary statistic $\hat\theta_{1:n}$ for the one-pass density estimator $\hat f$. \lin{A complete description of the above procedure is provided in Algorithm \suppref{alg:nonparametric-estimation}{1} in Section \suppref{sec:alg}{S5} of the supplement.}

\subsection{Tuning Parameter Selection}\label{subsec:tuning-parameter}
According to the pre-estimation strategy in Section \ref{subsec:increase-basis-functions}, the parameter $q$ is implicitly determined by a monotone function $S(\cdot)$. 
In addition, instead of starting with one basis function, we may start with $q_0>1$ basis functions to provide a better approximation to the regression function at the early stage of streaming. 
We therefore propose to set $S(j)=\lfloor (C_qj)^{1/h}\rfloor$, $\tau_j=\lfloor c_\circ S(j)\rfloor$ and $p=q=\min\{\mathfrak s_n/3,\max\{q_0,\lfloor S^{-1}(n)\rfloor\}\}$ for positive constants $h$, $C_q$, $c_\circ$ and  $q_0$, where $\mathfrak s_n$ represents the amount of available memory units at the timestamp $n$ for maintaining summary statistics; such choice of $p$ and $q$ respects the potential memory constraint $\mathfrak s_n$ ($\mathfrak s_n=\infty$ if there is no memory constraint). The constant $c_\circ$ controls the amount of pre-estimation and we recommend $c_\circ=1/2$ to provide sufficient pre-estimation for the basis coefficients. The numeric experiments in Section \suppref{sec:additional-simulation}{S7} of the supplement suggest that the performance is relatively robust to various choices of $C_q$ and $q_0$ around the recommended values $C_q=1/2$ and $q_0=5$; these values are also adopted in the simulation studies of Section \ref{subsec:simulation}.

\label{link:semi data-driven}
A fully data-driven approach to select the value of $h$ as well as $\rho$ in \eqref{eq:mhat-space-saving} is presented in Section \suppref{sec:additional-tuning}{S6} of the supplement. As the fully data-driven approach requires a large batch size $\batchsize$ to have a relatively robust choice of $\rho$,  we propose a semi data-driven strategy for the Fourier basis and the roughness penalty, in which case $\zeta = 4$ (defined in Assumption \ref{ass:W}), as follows; similar strategies can be derived for other bases and penalty matrices $W$. First, according to Corollary \ref{cor:L2-rate} we observe that $\rho\asympleq (q^{1-2\zeta}/n)^{1/2}\asymp n^{-(7h+1)/2}$ since $q\asymp S^{-1}(n)\asymp n^h$ based on our choice of the function $S(\cdot)$. We may then simply set $\rho = C_\rho n^{-(7h+1)/2}$ for a positive constant $C_\rho$, and determine both $h$ and $C_\rho$ based on the first $n_0$ observations by minimizing the cross-validation error
\begin{equation*}
   \mathrm{CV}(C_\rho,h)=\sum_{j=1}^J\sum_{i \in \mathscr{P}_j}\{Y_i-\hat m_{-j}^{C_\rho,h}(T_i)\}^2, 
\end{equation*}
where $\mathscr P_1,\ldots,\mathscr P_J$ for a positive integer $J$ form a partition of the first $n_0$ observations, and $\hat m_{-j}^{C_\rho,h}$ is  estimated with $q=\lfloor n^h/C_q\rfloor$ and $\rho=C_\rho n^{-(7h+1)/2}$ from the first $n_0$ observations except those in $\mathscr P_j$. This strategy ensures that $\rho$ decreases steadily in the long run while adapts to  the data noise level via the data-driven coefficient $C_\rho$. The simulation studies in Section \ref{subsec:simulation}  provide numeric illustrations of this simple and yet effective strategy. A numeric comparison presented in Section \suppref{sec:additional-tuning}{S6} of the supplement shows that the semi data-driven strategy is at least as effective as the fully data-driven approach, and thus is preferred in practice due to its simplicity and robustness.

\section{Minimax Convergence Rates}\label{sec:convergence rate}
This section is devoted to providing upper bounds on the uniform and $L^2$ convergence rates of the proposed space-saving one-pass estimator when only $O(q)$ memory units are available; matching lower bounds that consequently establish optimality of the proposed estimator are given in Section \ref{sec:space-lower-bound}.
Without loss of generality, we assume the batch size $B=1$.

\paragraph{Notation and convention.} For $k\in (0,\infty)$, let $L^k(\tdomain)$ be the collection of real-valued measurable functions $g$ such that $\|g\|_{L^k}:=\{\int_{\tdomain}|g(t)|^kdt\}^{1/k}<\infty$, and $L^\infty(\tdomain)$ the collection of real-valued measurable functions that are essentially bounded on $\tdomain$. Define $\|g\|_\infty=\sup_{t\in\tdomain}|g(t)|$ and note that $\|g\|_\infty=\|g\|_{L^\infty}$ when $g$ is continuous. For a square matrix $U$, $\lambda_{\min}(U)$ and $\lambda_{\max}(U)$ denote the minimum and maximum eigenvalues of $U$, respectively. We write $b_n\asympleq c_n$ or $b_n=O(c_n)$ (respectively, $b_n\asympgeq c_n$ or $b_n=\Omega(c_n)$) for two sequences $b_n$ and $c_n$ if there exists a constant $C>0$ such that $b_n\leq Cc_n$ (respectively, $b_n\geq Cc_n$) for all sufficiently large $n$, write $b_n\asymp c_n$ if and only if both $b_n\asympleq c_n$ and $b_n\asympgeq c_n$ hold, and write $b_n\ll c_n$ if $\lim_{n\rightarrow \infty} b_n/c_n=0$.

\subsection{Orthonormal Bases}\label{subsec:basis}

The properties of the estimator $\hat m$ depend on the characteristics of the adopted basis. We say a basis $\{\phi_j\}_{j=1}^\infty$ is a $(C,\alpha)$-basis for positive constants $C,\alpha$ if $\sup_{t\in\tdomain}\sum_{j=1}^q\phi_j^2(t) \leq C q^{\alpha}$, and a $(C,\alpha,\alpha_1)$-basis for an additional constant $\alpha_1$ if $\sup_{t\in\tdomain}\sum_{j=1}^q \{\phi_j^\prime(t)\}^2\leq Cq^{\alpha_1}$ also holds, for all $q\geq 1$. It turns out that the constant $\alpha$ plays a role in the convergence rate of the estimator $\hat m$. 

The approximation power of the basis also has impact on the convergence rate. For constants $\beta,\chi,M>0$, with $\nu$ being the largest integer less than $\beta$, let $\mathcal H(\beta,\chi,M)$ be the H\"older class 
\begin{align*}
    \mathcal H(\beta,\chi,M) =&\big\{m\in L^\infty(\tdomain): m \text{ is } \nu\text{-times continuously differentiable, }\|m\|_\infty \leq M \\ 
    & \qquad\qquad\qquad\,\, \text{and } |m^{(\nu)}(s)-m^{(\nu)}(t)|\leq \chi |s-t|^{\beta-\nu} \,\text{ for all } s,t\in\tdomain\big\}.
\end{align*}
 %$\mathcal H(\beta,\chi,M)$ denote the H\"older class of $\nu$-times continuously differentiable function $m$ such that $\|m\|_\infty \leq M$ and $|m^{(\nu)}(s)-m^{(\nu)}(t)|\leq \chi |s-t|^{\beta-\nu}$.
 Define  $\kappa_{\phi,2}^\ast(m,q)=\|m-P_qm\|_{L^2}$ and $\kappa_{\phi,2}(q)=\sup_{m\in\mathcal{H}(\beta,\chi,M)}\kappa_{\phi,2}^\ast(m,q)$ for an orthonormal basis $\{\phi_j\}_{j=1}^\infty$, where $P_qm:=\sum_{k=1}^q a_k\phi_k$ with $a_k=\int_{\tdomain} m(t)\phi_k(t)dt$. Similarly, define 
 $\kappa_{\phi,\infty}^\ast(m,q)=\|m- P_qm\|_\infty$ and  $\kappa_{\phi,\infty}(q)=\sup_{m\in\mathcal H(\beta,\chi,M)}\kappa_{\phi,\infty}^\ast(m,q)$. The quantities $\kappa_{\phi,2}(q)$ and $\kappa_{\phi,\infty}(q)$  characterize the approximation power of a given basis $\{\phi_j\}_{j=1}^q$ relative to $\mathcal H(\beta,\chi,M)$. \label{link:phi-revision}

Another feature of the basis $\{\phi_j\}_{j=1}^\infty$ relevant to the uniform convergence rate of $\hat m$ is the linear operator $\Pi_q^\phi$ defined by 
%\begin{equation*}
    $\Pi_q^\phi g = \phi^\top \sPhiPhi_q^{-1}\expect\{\phi(T)g(T)\}$
%\end{equation*}
 for any function $g$ with $\|g\|_\infty<\infty$, where $\phi=(\phi_1,\ldots,\phi_q)^\top$ and $\sPhiPhi_q$ is the $q\times q$ matrix formed by the elements $\int_{\tdomain} \phi_j(t)\phi_k(t)f(t)dt$ for $j,k=1,\ldots,q$. Define the $L^\infty$ operator norm of $\Pi_q^{\phi}$ by 
%\begin{equation*}
    $\|\Pi_q^{\phi}\|_\infty = \sup_{g:0<\|g\|_\infty<\infty}\frac{\|\Pi_q^{\phi} g\|_\infty}{\|g\|_\infty}.$
%\end{equation*}
Below are two examples of bases and their relevant features.%; additional bases are discussed in \cite{Chen2007} and \cite{belloni2015some}.

\begin{example}[Cohen–Daubechies–Vial wavelet series \citep{cohen1993wavelets}]\label{ex:wavelet-basis} The wavelet series,  $\tilde\phi_{J_0,0},\ldots,\tilde\phi_{J_0,2^{J_0}-1},\tilde\psi_{J_0,0},\ldots,\tilde\psi_{J-1,2^{J-1}-1}$,  described in Example 3.4 of \cite{belloni2015some}, for some $J_0$ satisfying $2^{J_0}\geq 2 \beta$, are  orthonormal  in $L^2([0,1])$.  Via suitable translation and scaling, they form an orthonormal $(C_\phi,1,\alpha_1)$-basis $\{\phi_j\}_{j=1}^\infty$ of $L^2(\tdomain)$ for some positive constants $C_\phi$ and $\alpha_1$. In addition, if $J=\log q$, then for a positive constant $C_{\beta,\chi}$ depending only on $\beta$ and $\chi$,  $\kappa_{\phi,2}(q)\leq C_{\beta,\chi} q^{-\beta}$ according to the discussions in Section 3.1of \cite{belloni2015some}. By Proposition 3.2 and Example 3.9 of \cite{belloni2015some}, $\|P_qm-m\|_\infty \leq (1+C_0)C_{\beta,\chi}q^{-\beta}$  for a constant $C_0>0$ and for all $m\in\mathcal H(\beta,\chi,M)$. This shows that $\kappa_{\phi,\infty}(q)\asympleq q^{-\beta}$. In addition, Theorem 5.1 of \cite{Chen2015} shows that $\|\Pi_q^{\phi}\|_\infty \asympleq 1$ for this basis when the density function $f$ is bounded away from zero and infinity.
\end{example}

\begin{example}[Fourier series]\label{ex:fourier-basis}
The Fourier basis functions $\phi_k$, defined by $\phi_1(t)=1$, $\phi_{2k}(t)=\cos(2k\pi t)$ and $\phi_{2k+1}(t)=\sin(2k\pi t)$ for $k\geq 1$, form an orthonormal basis of $L^2([0,1])$. Via suitable translation and scaling, they also form an orthonormal $(C_\phi,1,3)$-basis of $L^2(\tdomain)$ for a positive constant $C_\phi>0$.  In addition, for  a positive constant $C_{\beta,\chi}$ depending only on $\beta$ and $\chi$,  $\kappa_{\phi,2}^\ast(m,q)\leq  C_{\beta,\chi}q^{-\beta}$ and  $\kappa_{\phi,\infty}^\ast(m,q)\leq C_{\beta,\chi}(1+\log q)q^{-\beta}$ for all $\tdomain$-periodic functions $m$ in $\mathcal H(\beta,\chi,M)$, according to Eq. (5.8.4) of \cite{Canuto2006} and Example 3.7 of \cite{belloni2015some}. According to Example 3.7 of \cite{belloni2015some}, $\|\Pi_q^{\phi}\|_\infty\asympleq \log q$ when $f$ is the uniform distribution; see also Sections II.12 and III.13 of \cite{Zygmund2003}. In addition, the arguments in \cite{Zygmund2003} may be modified to accommodate densities that are bounded away from zero and infinity.
\end{example}

Compared to the wavelet basis, the Fourier basis has an extra $\log q$ factor in the rates of $\kappa_{\phi,\infty}(q)$ and  $\|\Pi_q^\phi\|_\infty$, which results in an additional $\log n$ factor in the uniform convergence rate of $\hat m$; see Corollary \ref{cor:fourier-uniform-rate}. However, due to its numerical stability, it is widely used in practice. Moreover, Fourier basis can be enhanced by a so-called Fourier extension technique to accommodate non-periodic functions; we provide a brief description of this technique in Section \suppref{sec:fourier-extension}{S11} of the supplement and also demonstrate its numerical performance in Section \ref{sec:numerical}.

\subsection{Assumptions}
We require some assumptions to study the theoretical properties of the space-saving one-pass estimator $\hat m$. The first one is on the design density of the predictor $T$, the second is on the smoothness of the regression function, and the third is on the noise variable.

\begin{assumption}[Design density]\label{ass:time-density}
The observations $T_{i}$ are independently and identically distributed (i.i.d.) according to a probability density function $f$ with $0<C_{f,1}\leq \inf_{t\in\tdomain} f(t)\leq \sup_{t\in\tdomain} f(t)\leq C_{f,2}<\infty$ for two positive universal constants $C_{f,1}$ and $C_{f,2}$.
\end{assumption}

\begin{assumption}[Regression regularity]\label{ass:regression-function}
The regression function $m$ belongs to $\mathcal H(\beta,\chi,M)$.
\end{assumption}

\label{page:R2:1}\begin{assumption}[Noise]\label{ass:noise} For universal constants $\xi>0$ and $C_\varepsilon>0$, 
the measurement errors $\varepsilon_i$ are {independent}, 
centered with $\expect{|\varepsilon_i|^{2+\xi}}\leq C_\varepsilon$, and independent of $T_1,T_2,\ldots$.
\end{assumption}

While Theorem \ref{thm:uniform-rate} holds for all $\xi>0$ in the above assumption, to achieve the optimal uniform convergence rate, we require $\xi\geq 1/\beta$; see Theorem \ref{thm:uniform-rate} and Corollary \ref{cor:uniform-rate} for details. In contrast, $\xi=0$ is sufficient for the optimal $L^2$ convergence rate; see Theorem \ref{thm:L2-convergence-rate} and Corollary \ref{cor:L2-rate}.

When the density function $f$ is unknown, the properties of the estimator $\hat m$ depend on the quality of the estimator $\hat f$, which in turn depends on the following regularity condition on $f$.
\begin{assumption}[Density regularity]\label{ass:density-regularity}
The density function $f$ belongs to the class $\mathcal D(\gamma,\chi,M)=\{f\in\mathcal H(\gamma,\chi,M):f\text{ is a probability density on }\tdomain\}$. 
\end{assumption}

\begin{remark} 
The assumed same constants $\chi$ and $M$ in  Assumptions \ref{ass:regression-function} and \ref{ass:density-regularity} are not essential, as the rate of convergence is influenced only by $\beta$ and $\gamma$. 
\end{remark}
%\vspace{5mm}

In the following condition we assume the penalty matrix $W$ to possess a non-negative spectrum that may grow with $q$. The simplest example of $W$ satisfying this condition is the identity matrix $I_q$. The roughness penalty matrix of the elements $\int_{\tdomain}\phi_j^{\prime\prime}(t)\phi_k^{\prime\prime}(t)dt$ also satisfies this condition, provided that the basis $\{\phi_j\}_{j=1}^\infty$ has the property $\sup_{t\in\tdomain}\sum_{j=1}^q \{\phi_j^{\prime\prime}(t)\}^2 \leq c_1 q^{c_2}$ for some constants $c_1$ and $c_2$;  the aforementioned wavelet and Fourier bases have this property. 
\begin{assumption}[Penalty matrix]\label{ass:W} For some constants $\zeta\geq 0$ and $C_{W,1}>0$, the matrix $W$ satisfies $0\leq \lambda_{\min}(W)\leq \lambda_{\max}(W)\leq C_{W,1}q^\zeta$.
\end{assumption}

\subsection{Convergence Rates}\label{subsec:rate}

Our first result concerns the uniform convergence rate of the one-pass density estimator $\hat f$. With $p\asymp (n/\log n)^{1/(2\gamma+1)}$, the uniform convergence rate of $\hat f$ matches the well known optimal rate of nonparametric density estimation \citep{stone1983optimal}. 
\begin{proposition}\label{prop:density-rate-Holder-class} If $\{\psi_j\}_{j=1}^\infty$ is a $(C_\psi,\alpha,\alpha_1)$-basis and {$p^{\alpha}\log(p+1)\asympleq n$}, then for any $\ell>0$, 
$$\sup_{f\in \mathcal D(\gamma,\chi,M)}\expect_f\|\hat f-f\|_\infty^\ell \leq c \{\sqrt{p^{\alpha}n^{-1}\log (n+1)} +  \kappa_{\psi,\infty}(p)\}^\ell,$$
where $\expect_f$ denotes the expectation with respect to the probability measure induced by the density function $f$, and $c$ is a constant depending only on $\ell,\gamma,\chi,M,C_\psi,\alpha,\alpha_1,c_\circ$.
\end{proposition}

\begin{corollary}
Under Assumption \ref{ass:density-regularity}, $\expect\|\hat f-f\|_\infty^2 \asympleq (n/\log n)^{-2\gamma/(2\gamma+1)}$ if one chooses $p\asymp (n/\log n)^{1/(2\gamma+1)}$ and $\{\psi_j\}_{j=1}^\infty$ is the wavelet basis in Example \ref{ex:wavelet-basis}.
\end{corollary}

 Let $\mathcal F=\mathcal F(\beta,\gamma,\chi,M)$ be the class of joint distributions on $(T_1,Y_1),\ldots,(T_n,Y_n)$  satisfying Assumptions \ref{ass:time-density}--\ref{ass:density-regularity}. When  $O(p+q)$ memory units are available so that we can use $O(p)$ basis functions to estimate the density function $f$ and $O(q)$ basis functions to estimate the regression function $m$, the theorem below establishes the uniform convergence rate of the space-saving one-pass estimator $\hat m$ uniformly over the class $\mathcal F$ under the memory constraint. \lin{In the following, recall that the concept of $(C,\alpha,\alpha_1)$-basis and the map $\kappa_{\phi,\infty}$ are introduced in Section \ref{subsec:basis}, and the constants $\xi$ and $\zeta$ are respectively defined in Assumption \ref{ass:noise} and Assumption \ref{ass:W}.}

%%%%%%%%%%%%%%%%%%%%%%%%%%%%%%%%%%%%%%%%%%
\begin{theorem}\label{thm:uniform-rate} Suppose both $\{\phi_j\}_{j=1}^\infty$ and $\{\psi_j\}_{j=1}^\infty$ are a $(C,\alpha,\alpha_1)$-basis. 
Assume  $p^{2\alpha}\asympleq n/\log n$,  {$q^\alpha\asympleq \min\{({n}/{\log n})^{\xi/(2+\xi)},(n/\log n)^{1/2}\}$}, $\kappa_{\phi,\infty}(q)\asympleq 1$, $\kappa_{\psi,\infty}(p)\ll 1$ and   
{$\rho^2 \asympleq (\log n)/(q^{2\zeta} n)$}.
Let {$r_n^2={p^{\alpha}(\log n)/n}$}, {$R_n^2={q^{\alpha}(\log n)/n}$} and $x_n=q^{\alpha/2}r_n^2+q^{\alpha/2}r_n \kappa_{\phi,\infty}(q)+q^{\alpha/2}\kappa_{\psi,\infty}(p)$.
Under Assumptions \ref{ass:time-density}--\ref{ass:W},  one has
$$\lim_{A\rightarrow\infty}\underset{n\rightarrow\infty}{\lim\sup}\sup_{F\in\mathcal{F}}\prob_F\big(\|\hat m- m\|_\infty \geq A\{r_n+R_n + (1+\|\Pi_q^{\phi}\|_\infty)\kappa_{\phi,\infty}(q)+x_n\}\big)=0,$$
where $F\in\mathcal F$ is a joint distribution on $(T_1,Y_1),\ldots,(T_n,Y_n)$, and $\prob_F$ denotes the probability measure induced by the distribution $F$.
\end{theorem}
%\begin{remark} \label{link:uniform-phase-transition}

\lin{In the above theorem, when $f$ is known, the terms $r_n$ and $x_n$ can be dropped along with Assumption \ref{ass:density-regularity} and  the conditions $q^{\alpha}\asympleq (n/\log n)^{1/2}$, $\kappa_{\psi,\infty}(p)\ll 1$ and $p^{2\alpha}\asympleq n/\log n$. When $f$ is unknown but sufficiently smooth, $r_n$ and $x_n$ are negligible relative to other terms.  The term $R_n$ is related to  estimation errors while $(1+\|\Pi_q^{\phi}\|_\infty)\kappa_{\phi,\infty}(q)$ stems from  approximation errors. When $q$ is sufficiently large, $R_n$ dominates the latter, and otherwise is dominated by the latter, with a phase transition at $q\asymp (n/\log n)^{1/(2\beta+1)}$ for the wavelet basis  as  $\|\Pi_q^{\phi}\|_\infty\asympleq 1$ for this basis. Also, shown in the following corollary, the proposed one-pass estimator achieves the  optimal uniform convergence rate in \cite{stone1982optimal} with $q\asymp (n/\log n)^{1/(2\beta+1)}$ given sufficient memory.}

%\vspace{-1cm}

\begin{corollary}[Estimating $m$ by the wavelet basis]\label{cor:uniform-rate} Suppose Assumptions \ref{ass:time-density}--\ref{ass:W} hold. If $\beta \geq 1/2$,  {$\gamma\geq \beta+1/2$}, $\xi\geq 1/\beta$, and both $\{\phi_j\}_{j=1}^\infty$ and $\{\psi_j\}_{j=1}^\infty$ are the wavelet basis in Example \ref{ex:wavelet-basis}, then with $q\asymp p\asymp(n/\log n)^{1/(2\beta+1)}$ and any $\rho \asympleq \{(\log n)/n\}^{(2\beta+2\zeta+1)/(4\beta+2)}$, one has $\|\hat m-m\|_\infty = \Op\big((n/\log n)^{-\beta/(2\beta+1)}\big).$
\end{corollary}

%\vspace{-1cm}

\lin{In the above corollary, the condition $\gamma\geq \beta+1/2$ explicitly requires the density function $f$ to be smoother than the regression function $m$; this requirement is also implicitly included in Theorem \ref{thm:uniform-rate} where we assume  $\kappa_{\psi,\infty}(p)\ll 1$ but $\kappa_{\phi,\infty}(q)\asympleq 1$. Certain violation of this condition is observed to have limited numeric impact; see Section \suppref{sec:disc-density}{S9} in the supplement for a numeric demonstration.}
If the Fourier basis is adopted and $m$ is periodic, the resulting estimator $\hat m$ achieves the optimal rate up to an additional $(\log n)^{(\beta+1)/(2\beta+1)}$ factor.

%\vspace{-1cm}

\begin{corollary}[Estimating $m$ by the Fourier basis]\label{cor:fourier-uniform-rate} Suppose Assumption \ref{ass:time-density}--\ref{ass:W} hold. Assume $\beta\geq 1/2$,  {$\gamma\geq \beta+1/2$}, $\xi\geq 1/\beta$ and  $\{\psi_j\}_{j=1}^\infty$ is the wavelet basis in Example \ref{ex:wavelet-basis}. If  $\{\phi_j\}_{j=1}^\infty$ is the Fourier basis in Example \ref{ex:fourier-basis} and $m$ is periodic on $\tdomain$, then with $q\asymp p\asymp n^{1/(2\beta+1)}$ and any $\rho \asympleq (\log n)^{1/2}/n^{(2\beta+2\zeta+1)/(4\beta+2)}$, one has $\|\hat m-m\|_\infty = \Op\big((n/\log n)^{-\beta/(2\beta+1)}(\log n)^{(\beta+1)/(2\beta+1)}\big).$
\end{corollary}

%\vspace{-1cm}

The following theorem establishes the $L^2$ convergence rate of the proposed estimator and its corollary shows that the estimator enjoys the optimal $L^2$ convergence rate \citep{stone1982optimal} provided sufficient memory. \lin{Similarly to the uniform convergence, there is a phase transition  at $q\asymp n^{1/(2\beta+1)}$ for the wavelet basis: The dominant term in the $L^2$ convergence rate is the approximation error $\kappa_{\phi,2}(q)$ when $q\ll n^{1/(2\beta+1)}$ and is the estimation error $R_n$ when $q\gg n^{1/(2\beta+1)}$.}
\begin{theorem}\label{thm:L2-convergence-rate}
Suppose $\{\phi_j\}_{j=1}^\infty$ is a $(C_\phi,\alpha)$-basis, $\{\psi_j\}_{j=1}^\infty$ is a $(C_\psi,\alpha,\alpha_1)$-basis, and Assumptions \ref{ass:time-density}--\ref{ass:W} hold, where $\xi=0$ in Assumption \ref{ass:noise} is allowed. 
Assume $p^{\alpha}n^{-1}\log(n)\ll 1$, $\kappa_{\psi,\infty}(p)\ll 1$, $\kappa_{\phi,\infty}(q)\asympleq 1$, $\rho n^{1/2}q^{\zeta-\alpha/2}\asympleq 1$, and {$q^\alpha\asympleq n$}. 
With $r_n^2=p^\alpha(\log n)/n$ and $R_n^2=q^\alpha/n$, one has 
$$\lim_{A\rightarrow\infty}\underset{n\rightarrow\infty}{\lim\sup}\sup_{F\in \mathcal F}\prob_F\big(\|\hat{m}-m\|_{L^2} \geq A\{r_n + R_n+\kappa_{\phi,2}(q)+\kappa_{\psi,\infty}(p)\}\big)=0.$$
\end{theorem}

%\vspace{-1cm}

\begin{corollary}\label{cor:L2-rate}
Under the assumptions of Theorem \ref{thm:L2-convergence-rate}, if $\gamma>\beta>0$ and both $\{\phi_j\}_{j=1}^\infty$ and $\{\psi_j\}_{j=1}^\infty$ are the wavelet basis in Example \ref{ex:wavelet-basis}, then with $p\asymp (n/\log n)^{1/(2\gamma+1)}$, $q\asymp n^{1/(2\beta+1)}$ and $\rho\asympleq (q^{1-2\zeta}/n)^{1/2}$, one has $\|\hat m-m\|_{L^2}=\Op(n^{-\beta/(2\beta+1)})$. If $\{\phi_j\}_{j=1}^\infty$ is the Fourier basis in Example \ref{ex:fourier-basis}, then the same convergence rate holds under the same conditions when $m$ is periodic on $\tdomain$.
\end{corollary}

%\vspace{-1cm}

\lin{The above theorems and their corollaries show that, when there are $\mathfrak s_n$ available memory units at the time $n$, if $\mathfrak s_n \to\infty$ and thus $p,q\rightarrow \infty$ is allowed as $n \to\infty$, then the estimate $\hat m$ converges to $m$ at certain rate. Otherwise, if there exists a constant $\mathfrak s_0$ such that $\mathfrak s_n \leq \mathfrak s_0<\infty$ for all $n$, then there is always a gap independent of $n$ between $\hat m$ and $m$.
As one of the main contributions, in the next section we show that this observation applies to all one-pass estimators that are able to produce an estimate of $m(t)$ simultaneously for all $t\in\tdomain$ at any time of the streaming process.}

\section{Space Complexity and Minimaxity}\label{sec:space-lower-bound}
While Corollaries \ref{cor:uniform-rate} and \ref{cor:L2-rate} show that the proposed estimator achieves the usual optimal rate when there is sufficient memory, a further important question concerns the best possible convergence rate \emph{under memory constraint} and whether the proposed estimator achieves this rate. When the amount of available memory units is capped by $b_n^{-1/\beta}$ for a converging sequence $b_n\asympleq (n/\log n)^{-\beta/(2\beta+1)}$, according to Theorem \ref{thm:uniform-rate}, by using a wavelet basis and by choosing $p\asymp q\asymp b_n^{-1/\beta}$, the proposed space-saving one-pass estimator satisfies the memory constraint and converges uniformly at the rate $\Op(b_n)$. 
In this section, we show that $b_n$ is the minimax rate that any one-pass estimator can achieve given only $O(b_n^{-1/\beta})$ memory units for estimating $m\in \mathcal H(\beta,\chi,M)$. Consequently, {the proposed one-pass estimator is statistically optimal under memory constraint}, and conversely,  has asymptotically minimal memory footprints among all one-pass estimators with the convergence rate $\Op(b_n)$.

\subsection{Space Complexity of Statistical Estimation}
Among all one-pass algorithms that solve an estimation problem with the same quality, those requiring a low level of memory footprints are generally preferred. For data analysts, given a specific analysis task, such as the one-pass nonparametric regression studied in this paper, of interest are the following dual questions, namely, whether it is possible to further reduce memory footprints without compromising statistical accuracy, and what the best statistical accuracy is given the amount of available memory. To answer these questions, one needs to find out the (asymptotically) tight lower bound on memory space that is required to achieve certain level of statistical accuracy. We say a statistical algorithm is optimal in memory usage if the amount of memory it requires matches the lower bound. While optimality in statistical accuracy, such as the minimax convergence rate of an estimation problem, has received considerable attention, few research has been devoted to optimality in  memory usage. Below we start to lay out a general framework for studying the interplay between statistical efficiency and memory usage.

Intuitively, the amount of demanded memory depends on the required estimation quality. For example, in the extreme case that we have no requirement on the estimation quality, then a constant memory space is sufficient, since for example the silly estimator $\breve m_n(t)\equiv Y_n$  requires only a fixed amount of memory units. Therefore, to formulate a sensible concept of space lower bound for an estimation problem, it is necessary to factor in a requirement on the estimation quality. Roughly speaking, a space lower bound may be defined for all estimators whose estimation error is (asymptotically) upper bounded by a quantity $b_n$ potentially varying with $n$.

To precisely formulate the above idea, we consider the following general one-pass estimation problem and later specialize it to one-pass nonparametric regression. Let $\mathcal P$ be a class of distributions on a sample space and $\vartheta:\mathcal P\rightarrow \Xi$ a function defined on $\mathcal P$ and taking values in a (potentially infinite-dimensional) parameter space $\Xi$. Given a data stream drawn from an unknown source distribution $P\in\mathcal P$, the goal is to estimate $\vartheta(P)$. 
For an estimator $\hat\vartheta$,  we write  $\hat\vartheta_n$ to denote its instance computed from the first $n$ observations of the data stream. With respect to a loss function $\ell:\Xi\times\Xi\rightarrow\real_{\geq 0}$, the estimation quality of $\hat\vartheta_n$ is quantified by the risk  $\expect_P\ell(\hat\vartheta_n,\vartheta(P))$, where $\expect_P$ denotes the expectation with respect to the distribution $P$. For a sequence $b=\{b_n\}_{n=1}^\infty$ of positive constants $b_n$, we define $\mathcal R_r(b,\mathcal P)=\{\hat\vartheta: \sup_{n\geq r}\sup_{P\in\mathcal P}b_n^{-1}\expect_P\ell(\hat\vartheta_n,\vartheta(P))\leq 1\}$ for any fixed $r\geq 1$. Roughly speaking, the maximum risk of any estimator in $\mathcal R_r(b,\mathcal P)$  converges at least at the rate $b_n$. In statistics, convergence rates in probability are also commonly used to quantify estimation quality. To accommodate this practice, we consider the class  $\mathcal E_r(b,\mathcal P,\delta)$ of estimators such that, if $\hat\vartheta\in \mathcal E_r(b,\mathcal P,\delta)$, then $\sup_{n\geq r}\sup_{P\in\mathcal P}\prob_P\{\ell(\hat\vartheta_n,\vartheta(P)) \geq b_n\} \leq \delta$. Intuitively, an estimator in the class  $\mathcal E_r(b,\mathcal P,\delta)$  converges at the rate $b_n$ with high probability. 
By Markov's inequality,  for any $\delta\in (0,1)$, it is seen that 
\begin{equation}\label{eq:R-E}
\mathcal R_r(b,\mathcal P)\subset \mathcal E_r(\delta^{-1}b,\mathcal P,\delta).   
\end{equation}

\begin{remark}
The constant $1$ in the definition of $\mathcal R_r(b,\mathcal P)$ is not essential; if it were some other positive constant $c$, then rescaling $b_n$ to $c^{-1}b_n$ would give rise to the same class $\mathcal R_r(b,\mathcal P)$.
\end{remark} 

Let $J_{\hat\vartheta}(z_1,\ldots,z_n)$ be the amount of memory units required to compute  $\hat\vartheta_n$ when the observed data stream is $(z_1,\ldots,z_n)$. Define $\mathcal J_n(\hat\vartheta,P)=\sup_{z_1,\ldots,z_n}J_{\hat\vartheta}(z_1,\ldots,z_n)$ with the supremum ranging over all possible sequences $(z_1,\ldots,z_n)$ sampled from the distribution $P$. The quantity $\sup_{P\in\mathcal P}\mathcal J_n(\hat\vartheta, P)$ then represents the amount of memory required for computing $\hat\vartheta_n$ in the worst case when the length of the data stream is $n$. The lower bound in the worst case over a class $\mathcal A$ of estimators  is denoted by $$\mathfrak{M}_n(\mathcal A, \mathcal P)=\inf_{\hat\vartheta\in\mathcal A}\sup_{P\in\mathcal P} \mathcal J_n(\hat\vartheta,P).$$ We say an estimator $\hat\vartheta\in\mathcal A$ is asymptotically optimal in memory usage within the class $\mathcal A$ if $\sup_{P\in\mathcal P}\mathcal J_n(\hat\vartheta,P)\asymp \mathfrak{M}_n(\mathcal A, \mathcal P)$. In analogy to  asymptotic minimaxity in statistics, such an estimator $\hat\vartheta$ may be also called an asymptotically minimax estimator in memory usage. Of interest are the lower bounds $\mathfrak{M}_n(\mathcal R_r(b,\mathcal P),\mathcal P)$ and $\mathfrak{M}_n(\mathcal E_r(b,\mathcal P,\delta),\mathcal P)$ for any fixed $r\geq 1$. In light of \eqref{eq:R-E}, it is straightforward to observe the following inequality, and therefore it is sufficient to consider the class $\mathcal E_r(b,\mathcal P,\delta)$  in the sequel.
\begin{proposition}
For any $\delta\in(0,1)$, one has
$\mathfrak{M}_n(\mathcal R_r(b,\mathcal P))\geq\mathfrak{M}_n(\mathcal E_r(\delta^{-1}b,\mathcal P,\delta))$.
\end{proposition}

\begin{remark}Since $\mathcal E_1(b,\mathcal P,\delta)\subset \mathcal E_2(b,\mathcal P,\delta)\subset\cdots$, it seems natural to consider the class $\mathcal E_\infty(b,\mathcal P,\delta)=\bigcup_{r=1}^\infty \mathcal E_r(b,\mathcal P,\delta)$. However, this class may be  too large to have a nontrivial lower bound. For example, let $\hat\vartheta^j$, $j=1,2,\ldots$, be a sequence of estimators such that each $\hat\vartheta^j$ requires $\lfloor q_n/j\rfloor+1$ memory units and $\lim\sup_n \sup_{P\in\mathcal P}\prob_P\{\ell(\hat\vartheta^j,\vartheta(P))\geq b_n\}\leq \delta/2$ for any $j\geq 1$, where $q_n$ is  a positive integer potentially growing with $n$; Theorem \ref{thm:uniform-rate} implies existence of such a sequence in nonparametric regression under some conditions. Since  $ \sup_{P\in\mathcal P}\mathcal J_n(\hat \vartheta^j, P)=\lfloor q_n/j\rfloor+1$ and $\hat\vartheta^j\in \mathcal E_{r_j}(b,\mathcal P,\delta)\subset \mathcal E_\infty(b,\mathcal P,\delta)$ for some $r_j<\infty$ and for each $j\geq 1$, we have the trivial lower bound $\mathfrak{M}_n(\mathcal E_\infty(b,\mathcal P,\delta))= 1$.
\end{remark}
%\vspace{5mm}

\subsection{Interplay between  Statistical Efficiency and Memory Constraint in One-pass Nonparametric Regression}\label{subsec:space-np-reg}

We now apply the general framework developed in the previous subsection to the one-pass nonparametric regression \eqref{eq:npreg} by taking $\Xi=\mathcal H(\beta,\chi,M)$, $\mathcal P=\mathcal F=\mathcal F(\beta,\gamma,\chi,M)$ of joint distributions of $(T_1,Y_1),\ldots,(T_n,Y_n)$, $\vartheta(P)=m(\cdot)= \expect_P(Y\mid T=\cdot)\in \mathcal H(\beta,\chi,M)$,   $\ell(g_1,g_2)=\|g_1-g_2\|_\infty$ for $g_1,g_2\in \mathcal H(\beta,\chi,M)$, and $\mathcal E_r(b,\mathcal F,\delta)$ to be the class of one-pass estimators $\breve m$ for $m$ such that $\sup_{n\geq r}\sup_{F\in\mathcal F}\prob_F\{\|\breve m-m\|_\infty\geq b_n\}\leq \delta$. 
\begin{remark}\label{rem:valid-estimator}
\lin{In this paper we consider only one-pass estimators $\breve m$ that are able to provide an estimate of $m(t)$ simultaneously for all $t\in\tdomain$ during the streaming process;  the OPLS estimator of \cite{yang2021online} is thus precluded according to Remark \ref{rem:diff-opls-proposed}.}
\end{remark}

As in \cite{Luo1993}, in the sequel we assume that one memory unit can hold any one of the real numbers. In reality, a computer represents a real number approximately by a floating-point data type up to a desired precision level, and in this case a memory unit corresponds to the chunk of memory required to hold a floating-point number. The following theorem establishes a lower bound on the memory footprints of any one-pass estimator for $m\in\mathcal H(\beta,\chi,M)$ with certain guarantee on estimation quality.
\begin{theorem}\label{thm:lower-bound-general-rate}
For a constant $\delta\in(0,1/2)$, if $\breve m$ is a one-pass estimator for the regression function  $m$ in \eqref{eq:npreg} such that 
\begin{equation}\label{eq:estimation-quality}
    {\sup_{n\geq r}}\sup_{m\in \mathcal H(\beta,\chi,M)}\prob_F(\|\breve m-m\|_\infty\geq b_n) \leq \delta
\end{equation}
for some $r\geq 1$, a sequence $b_n\asympleq 1$ and a joint distribution $F$ of $(T_1,\varepsilon_1),\ldots,(T_n,\varepsilon_n)$, then 
$\mathcal J_n(\breve m,F) \geq C b_n^{-1/\beta}$ for all $n\geq r$ and a constant $C>0$ not depending on $\breve m$.
%the amount of memory units required to compute $\breve m$ is $\Omega(b_n^{-1/\beta})$.
\end{theorem}
\begin{corollary}\label{cor:lower-bound-space}
$\mathfrak M_n(\mathcal E_r(Ab,\mathcal F,\delta),\mathcal F)\asympgeq b_n^{-1/\beta}$ for  any fixed $\delta\in(0,1/2)$, $A>0$ and $r\geq 1$.
\end{corollary}

The theorem stated below, directly following from Theorem \ref{thm:uniform-rate},
shows that the lower bound of Corollary \ref{cor:lower-bound-space} is tight for all sufficiently large $r$.
\begin{theorem}\label{thm:uniform-upper-bound}
In the regression problem \eqref{eq:npreg}, for $\beta\geq 1/2$, $\gamma\geq \beta+1/2$ and $\xi\geq 1/\beta$, for any $b_n\asympleq (n/\log n)^{-\beta/(2\beta+1)}$, there exists a one-pass estimator $\breve m$ for $m$ such that $\sup_{F\in\mathcal F}\mathcal J_n(\breve m,F)\asympleq b_n^{-1/\beta}$ and  \begin{equation}\label{eq:memory-upper}
\lim_{A\rightarrow\infty}\underset{n\rightarrow\infty}{\lim\sup}\sup_{F\in \mathcal F}\prob_F(\|\breve m-m\|_\infty \geq Ab_n) = 0.
\end{equation}
In particular, the space-saving one-pass estimator $\hat m$, estimated by using the first $p\asymp q\asymp b_n^{-1/\beta}$ wavelet basis functions, satisfies \eqref{eq:memory-upper} and $\sup_{F\in\mathcal F}\mathcal J_n(\hat m,F)\asympleq b_n^{-1/\beta}$.
%and only $O(b_n^{-1/\beta})$ memory units are required to compute $\hat m$.
\end{theorem}

The above results uncover some interesting interplay between space complexity and statistical efficiency in one-pass nonparametric regression. First, 
 there is a phase transition of the convergence rate at $\mathfrak S_n\asymp (n/\log n)^{1/(2\beta+1)}$ memory units: When the amount $\mathfrak s_n$ of available memory is asymptotically less than $\mathfrak S_n$, the uniform convergence rate of any one-pass estimator is restricted by the available memory and is capped at $\mathfrak s_n^{-\beta}$, and when there are at least $\mathfrak S_n$ memory units, the rate is restricted by the difficulty of nonparametric estimation  and is capped at the statistically optimal  rate $(n/\log n)^{-\beta/(2\beta+1)}$ of  \cite{stone1982optimal}. Consequently, in the regime $\mathfrak s_n\asympgeq \mathfrak S_n$, adding more memory cannot further improve the best possible rate.  On the other hand, if the amount of available memory is not scaled up  with the data volume, i.e., $\mathfrak s_n\leq \mathfrak s_0$ for some constant $\mathfrak s_0<\infty$,  then uniform convergence to $m$ is impossible for all one-pass estimators.

Second, any one-pass estimator achieving the optimal uniform convergence rate of \cite{stone1982optimal} requires at least $\mathfrak s_n\asympgeq (n/\log n)^{1/(2\beta+1)}$ memory units; this minimal memory requirement is achieved by the proposed estimator. Third, \lin{in the regime $1\ll \mathfrak s_n\ll \mathfrak S_n$, although the statistically optimal rate $(n/\log n)^{-\beta/(2\beta+1)}$ is not attainable, with $q\asymp \mathfrak s_n$, the proposed estimator converges to the underlying regression function at the best possible rate $\mathfrak s_n^{-\beta}\ll 1$  according to Theorem \ref{thm:uniform-upper-bound}. Consequently, the proposed  estimator $\hat m$ is asymptotically minimax in memory usage and  statistically optimal under memory constraint.}

For the $L^2$ loss function $\ell(g_1,g_2)=\|g_1-g_2\|_{L^2}$, we have the similar space lower bound stated in the following theorem. 
\begin{theorem}\label{thm:lower-bound-L2-rate}For any fixed $\delta\in(0,1/2)$, if  $\breve m$ is a one-pass estimator for the regression function $m$ in \eqref{eq:npreg} such that 
$$\sup_{n\geq r}\sup_{m\in \mathcal H(\beta,\chi,M)}\prob_F(\|\breve m-m\|_{L^2}\geq b_n) \leq \delta$$
for some $r \geq 1$, a sequence $b_n\asympleq 1$ and a joint distribution $F$ of $(T_1,\varepsilon_1),\ldots,(T_n,\varepsilon_n)$, then 
$\mathcal J_n(\breve m,F) \geq C b_n^{-1/\beta}$ for all $n\geq r$ and a constant $C>0$ not depending on $\breve m$.
\end{theorem}

The above lower bound is also asymptotically tight as it matches  the upper bound in the following theorem which follows directly from Theorem \ref{thm:L2-convergence-rate} with $\hat m$ and $\hat f$ estimated by using the first $p\asymp q\asymp b_n^{-1/\beta}$ wavelet basis functions. Note that the phase transition of the $L^2$ convergence rate  occurs at $\mathfrak S_n=n^{1/(2\beta+1)}$ memory units.
\begin{theorem}\label{thm:L2-upper-bound}
For the regression problem \eqref{eq:npreg}, for $\gamma>\beta>0$,  for any $b_n\asympleq n^{-\beta/(2\beta+1)}$, there exists a one-pass estimator $\breve m$ for $m$ such that $\sup_{F\in\mathcal F}\mathcal J_n(\breve m,F)\asympleq b_n^{-1/\beta}$ and
\begin{equation}\label{eq:memory-upper-L2}
\lim_{A\rightarrow\infty}\underset{n\rightarrow\infty}{\lim\sup}\sup_{F\in\mathcal F}\prob(\|\breve m-m\|_{L^2}\geq Ab_n) = 0.
\end{equation}
In particular, the space-saving one-pass estimator $\hat m$, estimated by using the first $p\asymp q\asymp b_n^{-1/\beta}$ wavelet basis functions, satisfies \eqref{eq:memory-upper-L2} and $\sup_{F\in\mathcal F}\mathcal J_n(\hat m,F)\asympleq b_n^{-1/\beta}$.
\end{theorem}

\begin{remark}\label{rem:opls-space-bound}
\lin{The OPLS estimator \citep{yang2021online} can produce a statistically optimal estimate of $m(t)$ for a predetermined finite set $\tdomain_\ast$ of $t$ by using $O(|\tdomain_\ast|)$ memory units, where $|\tdomain_\ast|$ denotes the number of elements in $\tdomain_\ast$. In addition, a numeric study in Section \suppref{sec:comp-time}{S10} of the supplement suggests that  OPLS  is efficient in answering queries about estimates of $m(t)$ for $t\in\tdomain_\ast$. 
These features make  OPLS  
%and thus may be 
suitable for the scenario that only predictions at $t\in\tdomain_\ast$ are needed. However, as per Remark \ref{rem:diff-opls-proposed}, a finite amount of memory is not sufficient for  OPLS  to estimate $m(t)$ simultaneously for all $t$. Therefore, in the application where an estimate of $m(t)$ for a random $t\in\tdomain$ may be requested at any time of the streaming process, our method may be preferred due to its optimality in memory usage and efficiency in computation time. 
A careful analysis of the spline-based method of \cite{yao2021spline}, which produces an estimate of $m(t)$ for all $t$, reveals that it requires asymptotically larger than $n^{1/(2\beta+1)}$ memory units to achieve the optimal convergence rate $n^{-\beta/(2\beta+1)}$ for estimating a $\beta$-times continuously differentiable function and thus is not optimal in memory usage according to Theorems \ref{thm:lower-bound-L2-rate} and \ref{thm:L2-upper-bound}.} 
\end{remark}

\subsection{Connection to Communication Complexity Theory}\label{subsec:ccct}

Our proof for Theorem \ref{thm:lower-bound-general-rate} exploits an intriguing connection between estimation quality and communication complexity. The theory of communication complexity, which is well established in computer science but less explored in statistics, concerns the amount of communication (measured by the number of bits or data units) required for two or more parties to solve a problem when the input to the problem is divided among the parties; see \cite{Rao2020} for an introduction. A problem particularly relevant to this paper is the following so-called one-way index problem involving two parties: One party, conventionally called Alice, holds a secret array $\omega$ of $k$ data units (e.g., an array of 0 and 1), and the other party, called Bob, holds a secret index $j$ and tries to find out the $j$th element $\omega_j$ of $\omega$ through a one-way communication protocol that only allows Alice to send messages to Bob. We say a randomized protocol solves the index problem with error probability $\delta\in(0,1/2)$ if the probability that Bob outputs a wrong value is at most $\delta$. For a given $\delta\in(0,1/2)$, the randomized one-way communication complexity of the index problem is the minimum worst-case (over all possible inputs and all possible realizations of the protocol) number of data units sent through any randomized one-way protocol that solves the index problem with an error probability bounded by $\delta$. A remarkable theorem in \cite{Kremer1999} asserts that such communication complexity is $\Omega(k)$ for any fixed $\delta\in(0,1/2)$.

The key idea to establish Theorem \ref{thm:lower-bound-general-rate} is that any estimator with the quality guarantee \eqref{eq:estimation-quality} can be used to construct a randomized protocol to solve the index problem. As in \cite{Luo1993}, we count the number of real numbers to be communicated in the communication complexity, i.e., each data unit is a real number. Suppose that Alice holds an element $\omega$ from  
 the $k$-dimensional pruned hypercube  $\mathbb H_k=\{(\omega_1,\ldots,\omega_k):\omega_j\in\{0,1\}\,\, \forall\,j=1,\ldots,k\}$ endowed with the Hamming distance $d(\omega,\nu)=\sum_{j=1}^k |\omega_j-\nu_j|$ for $\omega,\nu\in\mathbb H_k$. Below  we design a randomized one-way protocol for the index problem by using any estimator $\breve m$ satisfying \eqref{eq:estimation-quality}.

For a positive integer $k$ whose value is to be determined later, divide $\tdomain$ into $k$ intervals of the same length and let $t_1,\ldots,t_k$ be the centers of these intervals. Let $K$ be any sufficiently smooth function with support on $(-1/2,1/2)$ and satisfying $\sup_{t\in(-1/2,1/2)}|K(t)|= M/\chi$ and $K(0)>0$. For each $\omega\in \mathbb H_k$, we define $m_\omega(t)=\sum_{j=1}^k \omega_j c_K\chi k^{-\beta} K(k(t-t_j))$, where the positive constant $c_K$, depending only on $K$, is chosen to be sufficiently small so that $m_\omega\in\mathcal H(\beta,\chi,M)$. Then one can show that $\|m_\omega-m_\nu\|_\infty\geq cM k^{-\beta}$ with $c=c_K$ whenever $\omega\neq \nu$.

\begin{figure}
    \centering
    \begin{tikzpicture}[alicestyle/.style={draw,fill=red!30},bobstyle/.style={draw,fill=blue!30},orangestyle/.style={draw,fill=orange!30},scale=0.95,every node/.style={scale=0.95}]
	\def\aliceW{5.8}
	\def\aliceH{2.5}
	\def\bobW{4.2}
	\def\bobH{2.5}
	\def\boardW{7.2}
	\def\boardH{0.8}
	\def\D{2.1}
	\def\marg{0.02}
	
	\filldraw[alicestyle,rounded corners] (-\D-\aliceW,0) rectangle ++(\aliceW,-\aliceH);
	
	\node at (-\D-\aliceW/2,-0.7*\aliceH/5) {\small Alice: $\omega=(\omega_1,\ldots,\omega_k)$};
	\draw (-\D-\aliceW+\marg*\aliceW,-1.3*\aliceH/5) -- (-\D-\marg*\aliceW,-1.3*\aliceH/5);
	
	\node at (-\D-\aliceW/2,-2.2*\aliceH/5) {\small sample $(T_1,Y_1),\ldots,(T_n,Y_n)$};
	\node at (-\D-\aliceW/2,-3.2*\aliceH/5) {\small from the model  $Y_i=m_\omega(T_i)+\varepsilon_i$};
	\node at (-\D-\aliceW/2,-4.2*\aliceH/5) {\small and  produce the estimate $\breve m_\omega$};
	
	\draw[-stealth,thick] (-\D,-\aliceH/2) -> (\D,-\bobH/2);
	\node at (0,-0.9*\aliceH/6) {\small send the memory};
	\node at (0,-2*\aliceH/6) {\small footprints of $\breve m_\omega$};
	
	\draw[-stealth,thick] (\D+\bobW,-\aliceH/2) -> (\D+1.6*\bobW,-\aliceH/2);
	\node at (\D+1.3*\bobW,-2*\aliceH/6) {\small output $\breve \omega_j$};

	\filldraw[bobstyle,rounded corners] (\D,0) rectangle ++(\bobW,-\bobH);
	
	\node at (\D+\bobW/2,-0.7*\bobH/5) {\small Bob: $j$};
	\draw (\D+\marg*\bobW,-1.3*\bobH/5) -- (\D+\bobW-\marg*\bobW,-1.3*\bobH/5);
	
	\node at (\D+\bobW/2,-2.5*\bobH/5) {\small compute $\breve\omega=$};
	\node at (\D+\bobW/2,-3.5*\bobH/5) {\small  $\arg\min_{\nu\in\mathbb H_k} d(\breve m_\omega,m_\nu)$};
	
	\filldraw[orangestyle,rounded corners] (-\boardW/2,-1.7*\aliceH) rectangle ++(\boardW,\boardH);
	\node at(0,-1.7*\aliceH+\boardH/2) {\small Shared information: $\mathcal M_k=\{m_\omega:\omega\in\mathbb H_k\}$};

	\draw[-stealth,thick] (-\boardW/3,-1.7*\aliceH+\boardH) -> (-1.5*\boardW/3,-\aliceH);
	\draw[-stealth,thick] (\boardW/3,-1.7*\aliceH+\boardH) -> (1.5*\boardW/3,-\aliceH);
	
\end{tikzpicture}
    \caption{Illustration of the one-way protocol to solve the index problem.}
    \label{fig:diagram-of-communication}
\end{figure}
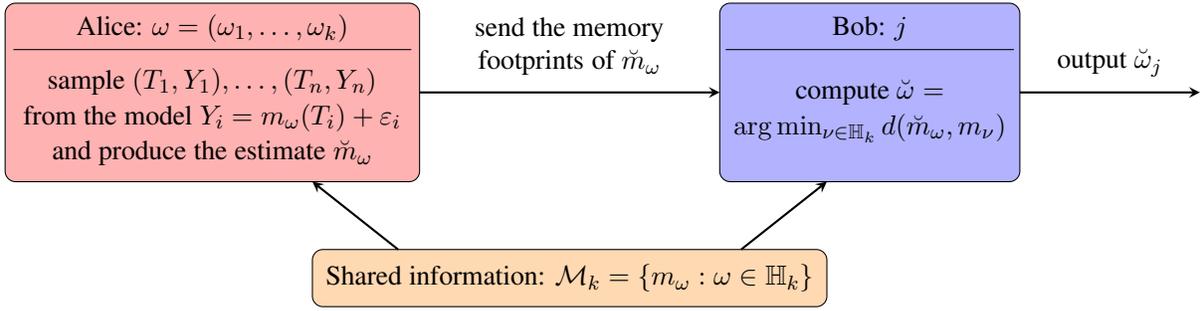

As a part of the protocol, the class $\mathcal M_k=\{m_\omega:\omega\in\mathbb H_k\}$ is known to both Alice and Bob. The protocol, illustrated in Figure \ref{fig:diagram-of-communication},  proceeds as follows.
First, Alice uses her secret $\omega$ to generate  independent observations $(T_1,Y_1),\ldots,(T_n,Y_n)$ according to the model $Y_i=m_\omega(T_i)+\varepsilon_i$ with $(T_1,\varepsilon_1),\ldots,(T_n,\varepsilon_n)$ following the joint distribution $F$ in Theorem \ref{thm:lower-bound-general-rate}. Based on these observations Alice uses an estimator $\breve m$ satisfying \eqref{eq:estimation-quality} to produce an estimate $\breve m_\omega$ and sends the memory footprints of $\breve m_\omega$ to Bob through a one-way communication channel. Upon receiving $\breve m_\omega$, Bob computes  $\breve\omega=\arg\min_{\nu\in\mathbb H_k} d(\breve m_\omega,m_\nu)$ and outputs $\breve\omega_j$ as an estimate of $\omega_j$.

We now show that, if $k$ is the largest integer that is not larger than $\{cM/(3b_n)\}^{1/\beta}$ so that   $3b_n\leq cMk^{-\beta}$, then this protocol correctly solves the index problem with error probability at most $\delta$ for all $n\geq r$. To see this, we first observe that when the event $\|\breve m_\omega - m_\omega\|_\infty < b_n$ occurs, for all   $\nu\neq\omega$, $\|\breve m_\omega - m_\nu\|_\infty \geq \|m_\omega-m_\nu\|_\infty-\|\breve m_\omega - m_\omega\|_\infty \geq 2b_n$, and consequently $\breve\omega=\omega$. By assumption \eqref{eq:estimation-quality}, the event $\|\breve m_\omega - m_\omega\|_\infty < b_n$ occurs with probability at last $1-\delta$, and thus Bob correctly outputs $\omega_j$ with probability at least $1-\delta$. %, or equivalently, the probability that Bob outputs a wrong value is at most $\delta$. 

Finally, according to \cite{Kremer1999}, the amount of data units transmitted, which is the amount of memory footprints of $\breve m_\omega$ according to the designed protocol, must be $\Omega(k)$, and Theorem \ref{thm:lower-bound-general-rate} follows by observing that the choice of $k$ implies $k\asymp b_n^{-1/\beta}$ and that all constants involved do not depend on any specific estimator $\breve m$ as long as it satisfies \eqref{eq:estimation-quality}. The proof for Theorem \ref{thm:lower-bound-L2-rate} is almost identical to that for Theorem \ref{thm:lower-bound-general-rate} and thus is omitted.

%%%%%%%%%%%%%%%%%%%%%%%%%%%%%%%%%%%%%%%%%%%%%%%%%%%%%%%%

\section{Numerical Studies}\label{sec:numerical}

\subsection{Simulation Studies}\label{subsec:simulation}
To illustrate the numerical performance of the proposed method, we conduct simulation studies with the following three settings for the regression function, namely,
\begin{description}
    \item[I.] $m_1(t)=\exp\{\sin(2\pi t)\}$, a smooth periodic function, 
    \item[II.] $m_2(t)=|t-0.4|$, which is a non-periodic function and non-differentiable at $t=0.4$ but smooth elsewhere, 
    \item[III.] $m_3(t)=\sum_{k=1}^\infty k^{-1.5}\phi_k(t)$ with $\phi_1(t)=1$, $\phi_{2k}(t)=\cos(2k\pi t)$ and $\phi_{2k+1}(t)=\sin(2k\pi t)$, a continuous periodic  function that is not twice differentiable everywhere.
\end{description}

In reality, data items often come or are processed in batches for computational efficiency. To mimic this practice, we set the batch size $B=100$. The observations $T_{i}$ of the predictor are independently sampled from the uniform distribution on $\tdomain=[0,1]$, while the measurement errors $\varepsilon_i$ are independently sampled from the centered Gaussian distribution {$N(0,\sigma^2)$} with $\sigma^2$ being determined by the signal-to-noise ratio $\mathrm{SNR}=2$ that is defined by $\mathrm{SNR}=\expect|m(T)|^2/\sigma^2$; non-uniform distributions and non-Gaussian noise are also investigated but not presented here, as they lead to similar results.
Estimation quality is quantified by root mean integrated squared error %(RMISE) 
%\begin{equation*}\label{eq:RMISE}
$\text{RMISE}=\sqrt{\frac{1}{N}\sum_{i=1}^N \int_\tdomain |m(t)-\hat{m}(t)|^2\,dt},$
%\end{equation*}
where $N=100$ is the number of Monte Carlo replicates in the studies. For the proposed method, we adopt the Fourier basis for both regression function  and density estimators, with the Fourier extension technique {(briefly described in Section \suppref{sec:fourier-extension}{S11} of the supplement)} to handle non-periodicity, where the extension margin is set to $0.1$.

We compare the proposed method with two recently developed one-pass nonparametric estimation methods, namely, the one-pass local smoothing (OPLS) estimator in \cite{yang2021online} and the dynamic penalized spline regression (DPSR) estimator in \cite{yao2021spline}.  \lin{For the proposed method, we determine the tuning parameters $q$ and $\rho$ by the semi data-driven strategy described in Section \ref{subsec:tuning-parameter}; note that $q$ grows with the sample size $n$.} 
\lin{For the OPLS method, the bandwidth is tuned by the method described in Section 3 of \cite{yang2021online} with the parameter $L=10$. For the DPSR method, cubic splines are adopted, where the parameter $\lambda_n$ is tuned by the generalized cross-validation method and the dynamic knots are updated according to the implementation proposed in Section 2.2 of \cite{yao2021spline} with the parameters $\alpha=1$ and $\nu=1/3$.} In addition, in the comparison we include the non-streaming estimator \eqref{eq:a-full} that is expected to yield a better performance since it has access to all historical data.

\begin{figure}
\centering
\subfloat[$m_1$]{
\includegraphics[scale=0.25]{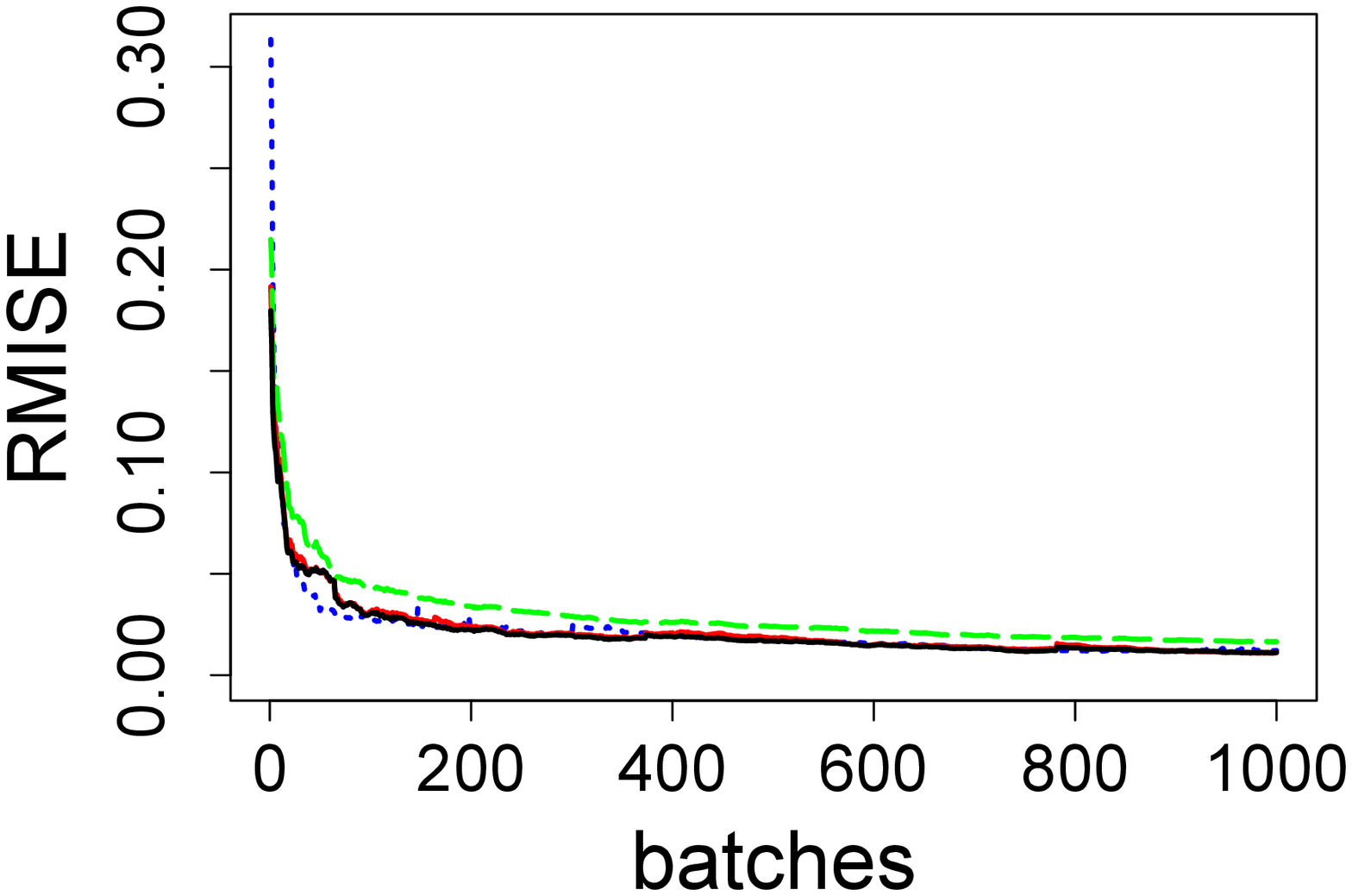}}\hspace{0.4cm}
\subfloat[$m_2$]{
\includegraphics[scale=0.25]{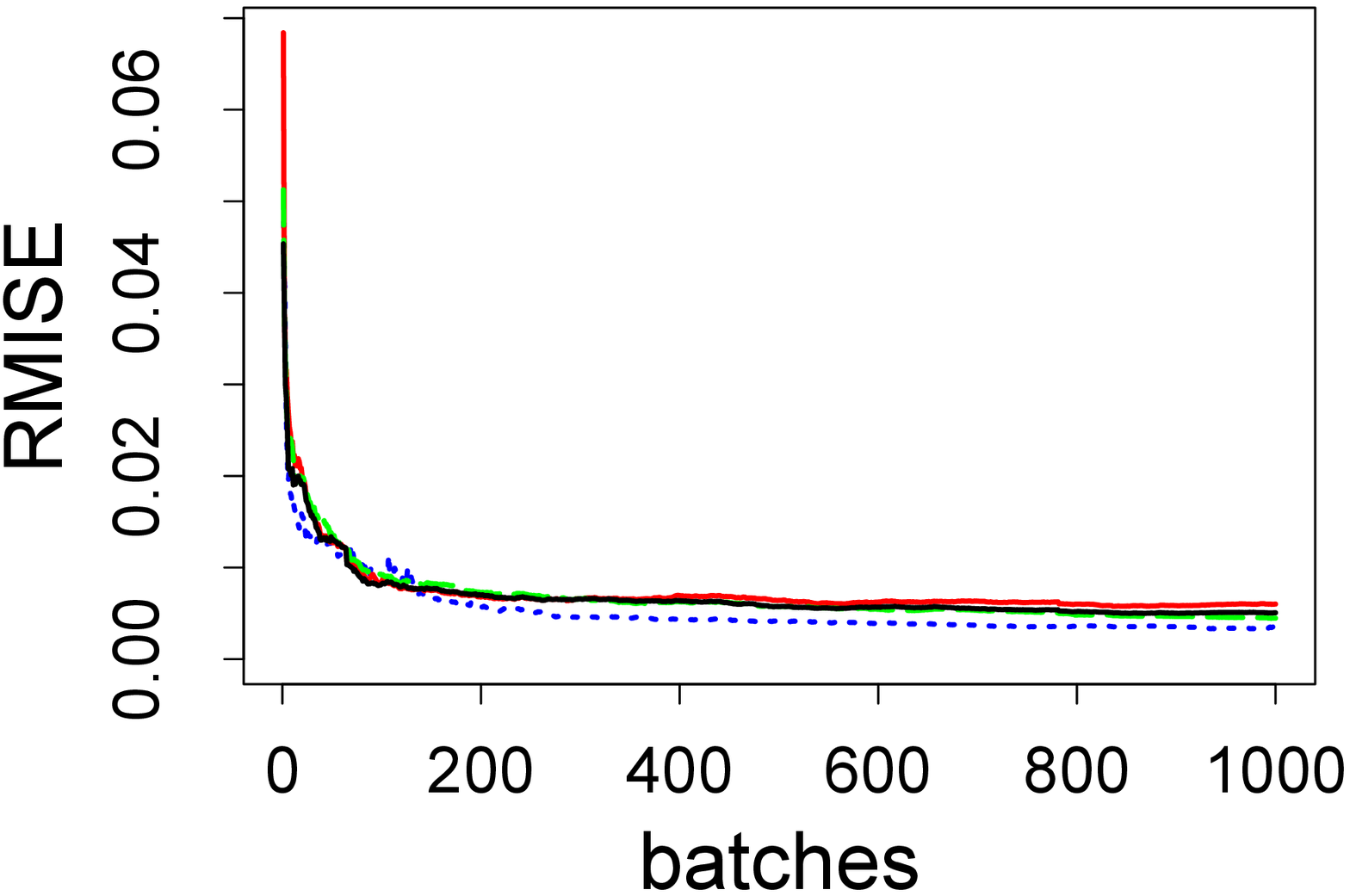}}\hspace{0.4cm}
\subfloat[$m_3$]{
\includegraphics[scale=0.25]{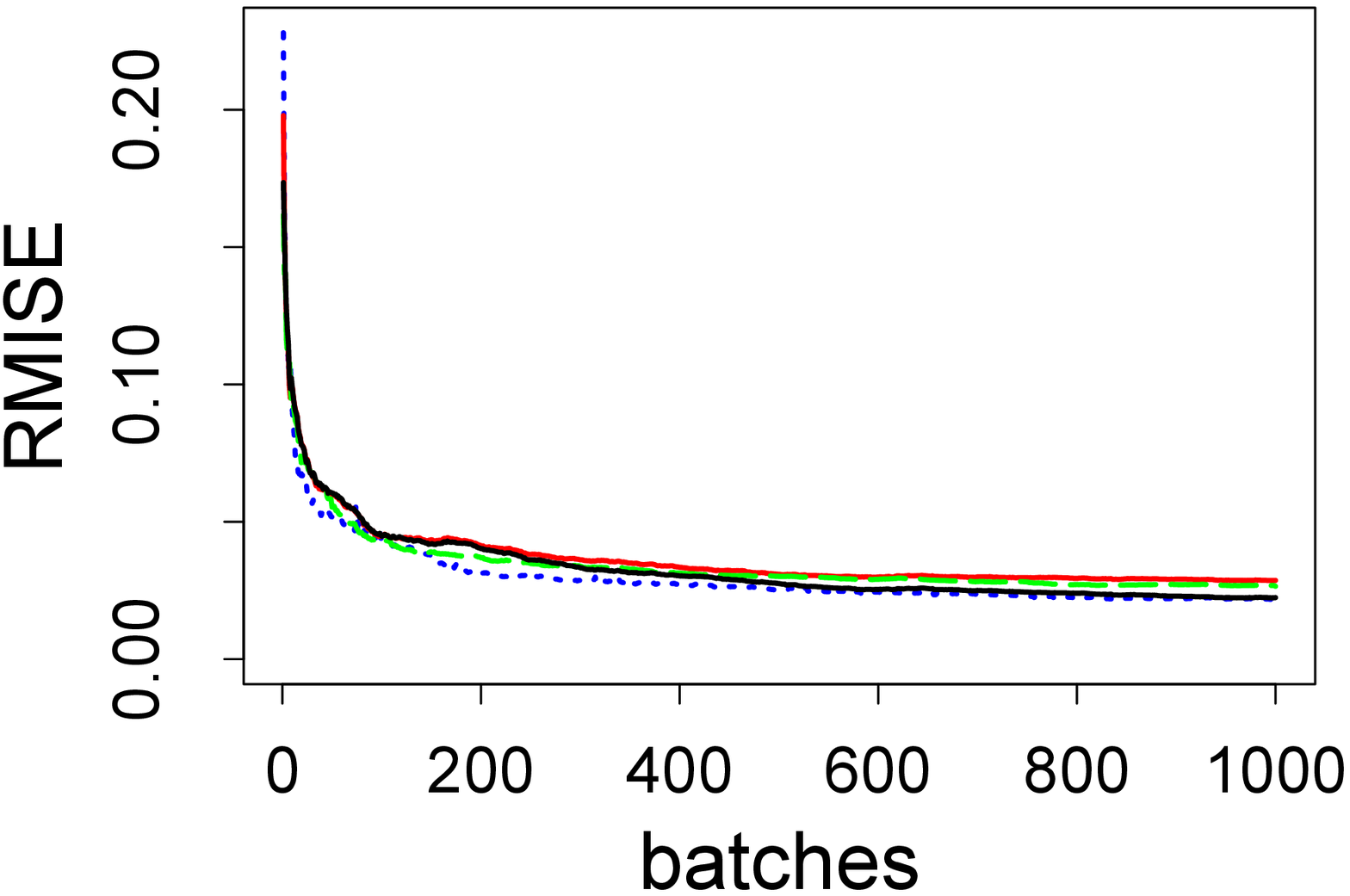}}
\caption{RMISE under different values of $n$ for the proposed method (red solid), the non-streaming method (black dashed-dotted),  OPLS  (green dashed) and  DPSR (blue dotted).}
\label{fig: RMISE-comparison-snr2}
\end{figure}

The results, presented in Figure \ref{fig: RMISE-comparison-snr2}, show that all estimators yield comparable performance in the long run. In particular, we observe that the proposed one-pass method is rather efficient relative to its non-streaming counterpart. The number of basis functions used in the proposed method for estimating smooth functions is considerably less than that of the DPSR method. For example, our estimator requires around 20 basis functions for estimating $m_2$ at the timestamp $n=100000$, while the DPSR method uses about 42 spline basis functions.   \lin{A comparison on computation time is presented in Section \suppref{sec:comp-time}{S10} of the supplement, which shows that the proposed estimator is also efficient in computation time.}  The convergence rates and phase transitions in Section \ref{subsec:rate} are numerically validated in Section \suppref{sec:numeric-comp-rates}{S8} of the supplement.

\subsection{A Case Study}\label{subsec:realdata}
%Streaming data are nowadays ubiquitous in many domains, such as network analytics,  business analytics and health/fitness monitoring.  
%As an example, 
Consider a fitness app used by a group of users to track and share their physical activities collected from wearable  devices. Relevant to the users, for example, are the contrast between their individual activity profiles and the average profile of certain cohort that the individual users belong to; users may utilize such information to adjust their activity levels. %, such as getting more exercise if their physical activity levels are much lower than the cohort average. 
In this scenario, data constantly generated by each device are not only sent to a server for central processing but also locally processed in the device, giving rise to two data streams, namely, the local stream processed in the device, and the global stream that collects observations from all devices and is processed in the server. Due to the massive data received by the server and the limited computation capacity of local devices, one-pass algorithms are preferred in both server and local devices.

For an illustration, we utilize the NHANES\footnote{https://wwwn.cdc.gov/nchs/nhanes} (National Health and Nutrition Examination Survey) 2005-2006 dataset to mimic data streams in a fitness app; although this dataset had been fully collected by 2006, it naturally formed data streams during the collecting process. In the study, each participant of age $6+$ years was asked to wear a physical activity monitor  during non-sleeping hours for seven days. For each minute, the monitor recorded the average of the physical intensity levels ({ranging from 0 to 32767}) within the minute. Consequently, a local data stream is observed for each participant, and these streams are then merged into a global stream. To further mimic the reality that observations (even recorded at the same time) from different users may arrive at the server in a random order potentially due to unpredictable network latency, we randomly interlace the local streams when merging them into the global stream.

\begin{figure}
\centering
\subfloat[$\lfloor n/3\rfloor$]{
\includegraphics[scale=0.22]{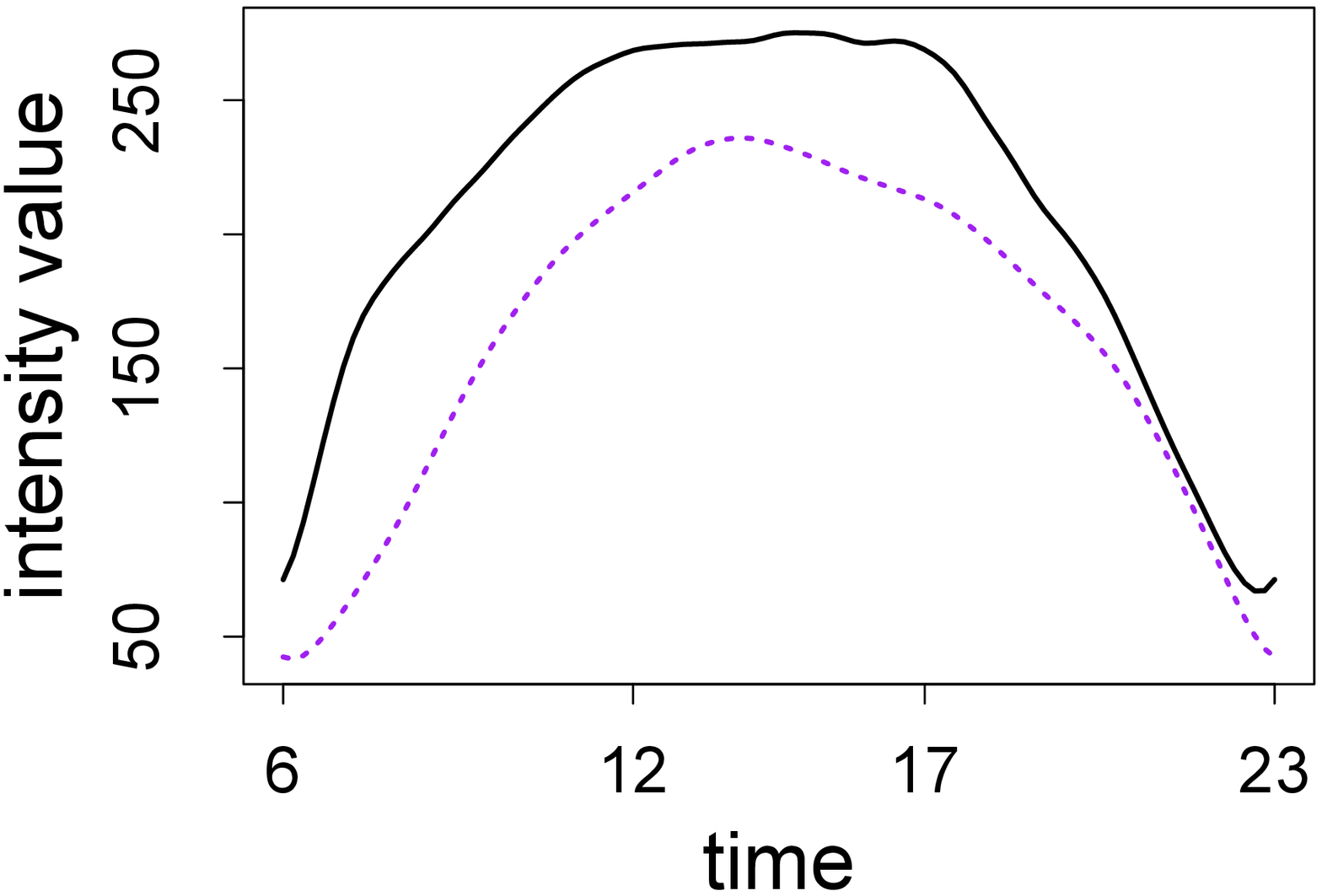}}\hspace{0.4cm}
\subfloat[$\lfloor 2n/3\rfloor$]{
\includegraphics[scale=0.22]{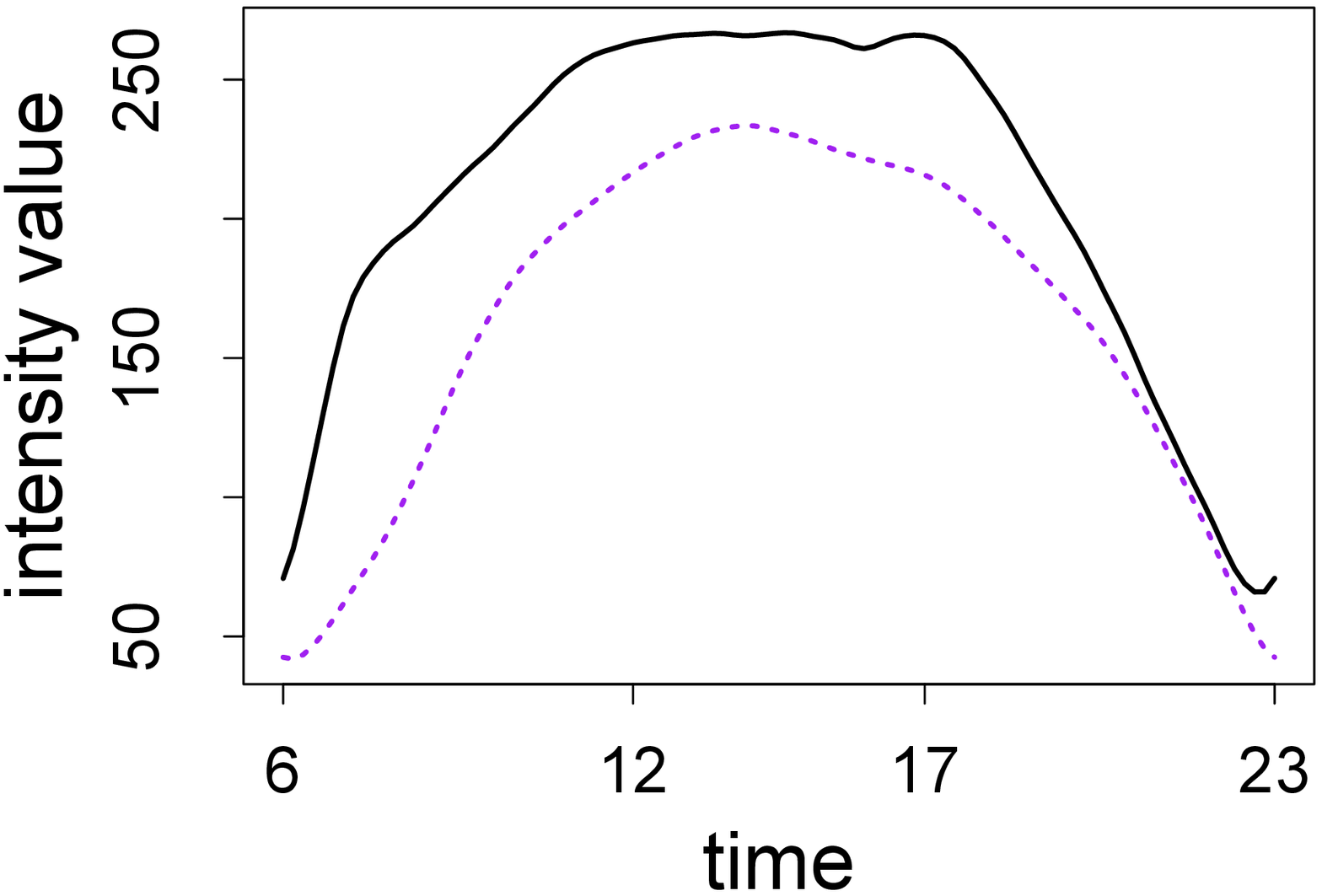}}\hspace{0.4cm}
\subfloat[$n$]{
\includegraphics[scale=0.22]{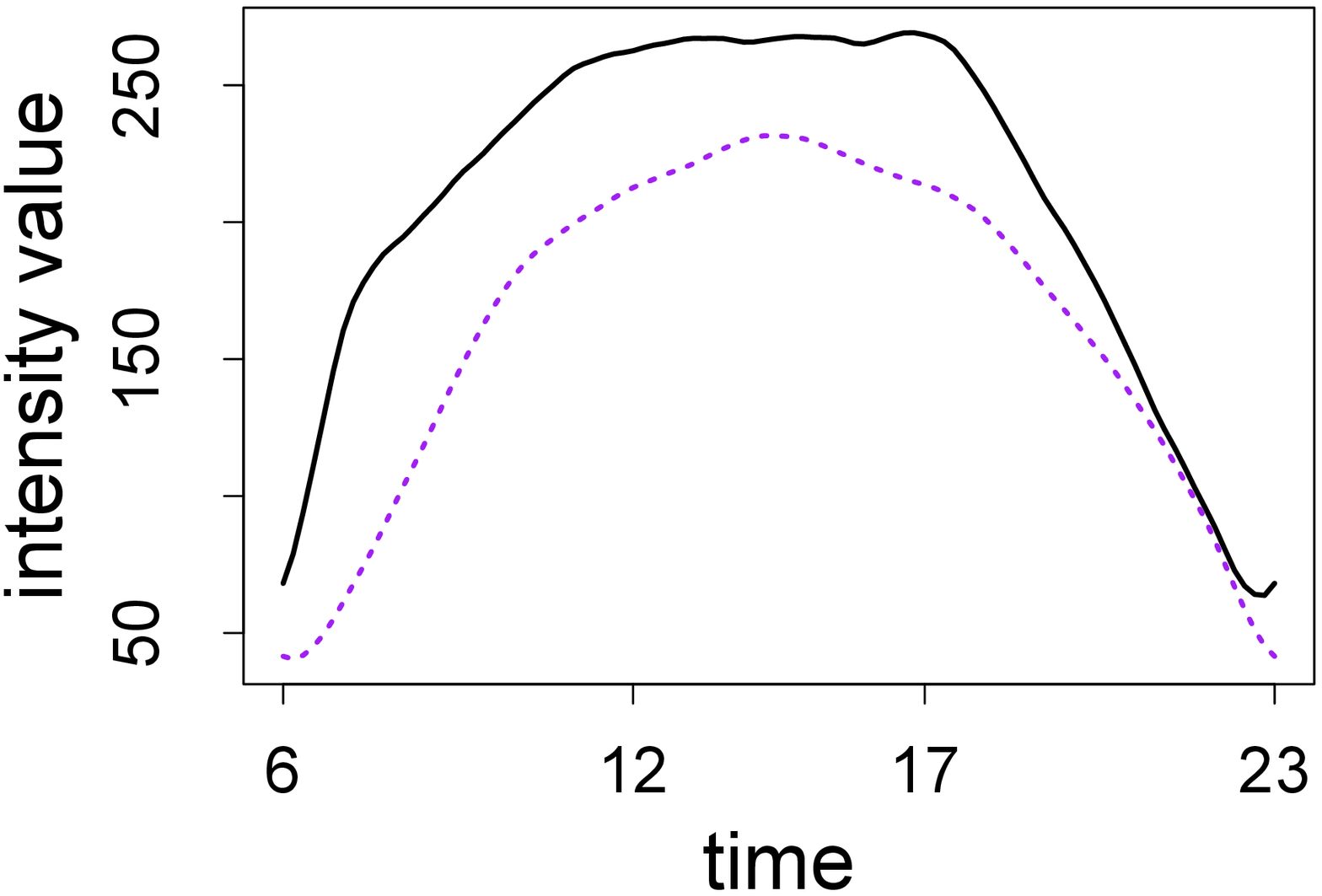}}
\caption{Snapshots of the dynamic weekday (black solid) and weekend (purple dotted) cohort physical activity profiles estimated by  the proposed method at various timestamps.}% $i=\lfloor n/3\rfloor,\lfloor 2n/3\rfloor,n$.}
\label{fig:real-data-population-comparison}
\end{figure}

We first focus on the cohort of adults of ages between 18 and 45 and their physical activities during the general non-sleeping hours 06:00--23:00. After excluding observations whose reliability is questionable according to the NHANES protocol, the global stream for this cohort contains $n=18330709$ observations $(T_i,Y_i)$, %from 2426 participants, 
where $Y_i$ is the intensity level recorded at the time $T_i\in [6,23]$, for $i=1,\ldots,n$. This stream is further divided into two substreams, one for the weekday and the other for the weekend, since human activities are often rather different in the weekday and weekend. Each substream is then fed into the proposed one-pass algorithm, where the value $\rho$ is determined by the semi data-driven strategy proposed in Section \ref{subsec:tuning-parameter}.

Shown in Figure \ref{fig:real-data-population-comparison} are the snapshots of the dynamic weekday and weekend cohort physical activity profiles at three different timestamps, namely, $i=\lfloor n/3\rfloor$, $\lfloor 2n/3\rfloor$ and $n$. Not surprisingly, the weekday cohort activity intensity level is uniformly higher than its weekend counterpart, suggesting that people tend to relax during the weekend.

\begin{figure}[t]
\centering
\subfloat[31212, weekday]{
\includegraphics[scale=0.2]{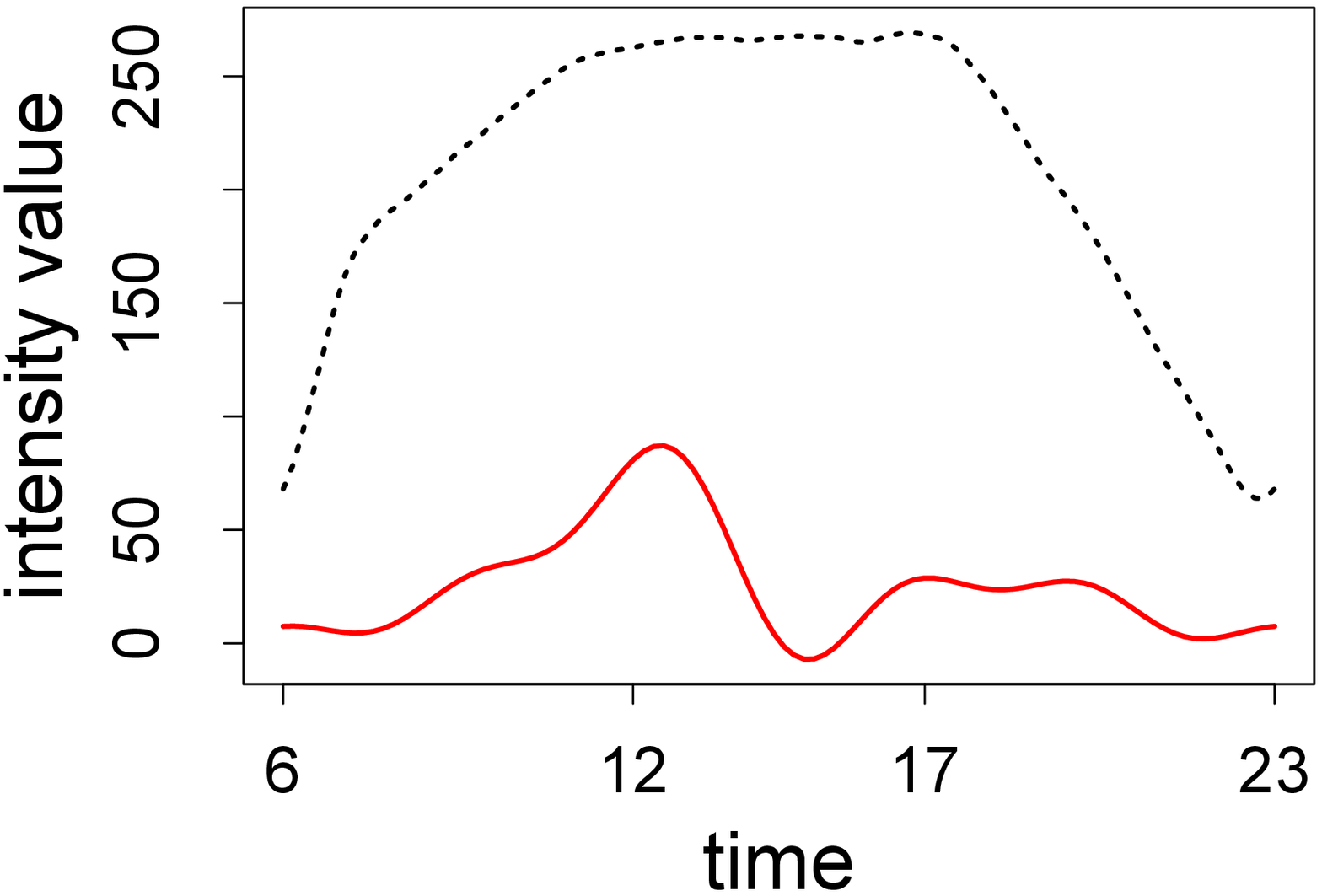}}\hspace{0.2cm}
\subfloat[31212, weekend]{
\includegraphics[scale=0.2]{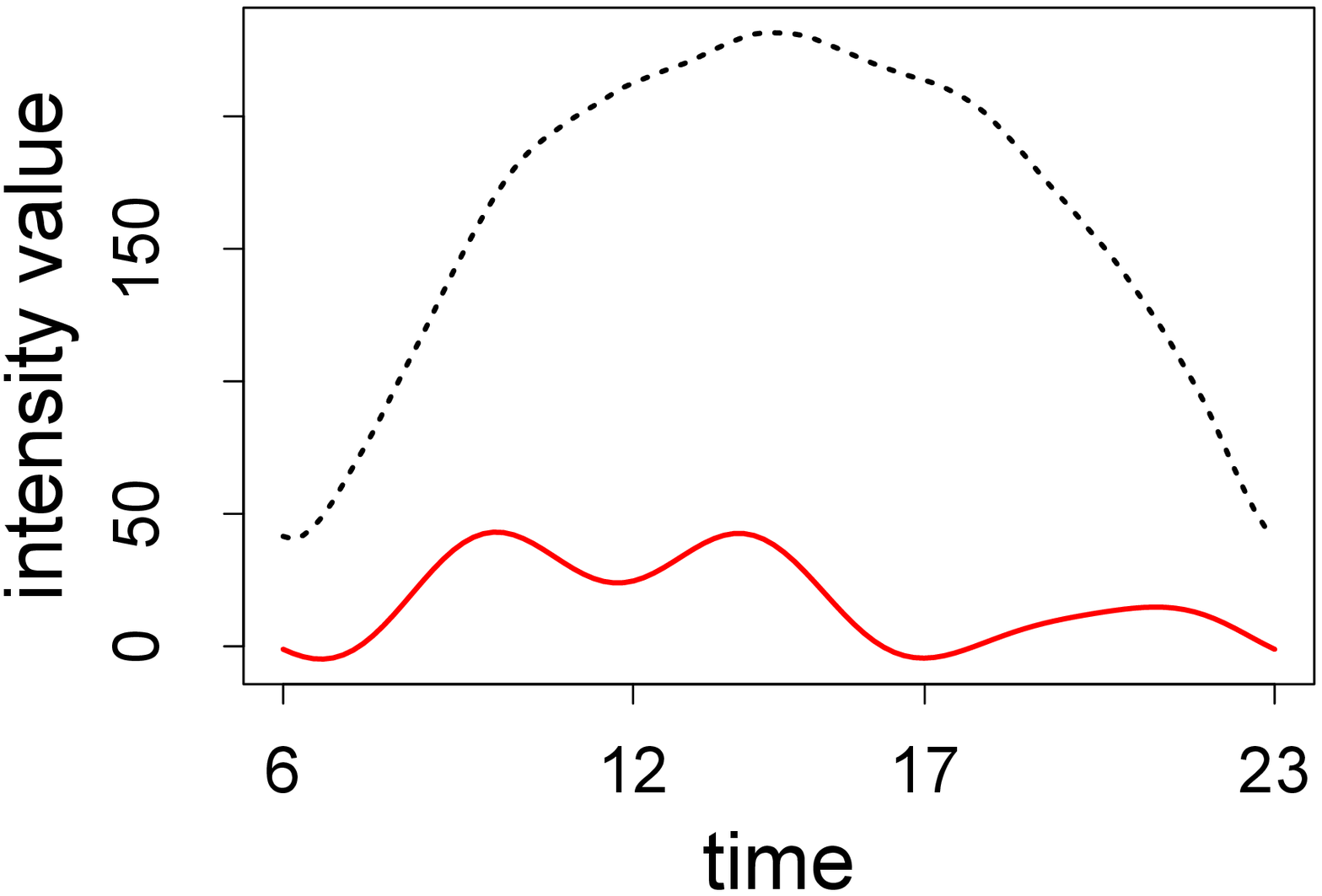}}\hspace{0.2cm}
\subfloat[33166, weekday]{
\includegraphics[scale=0.2]{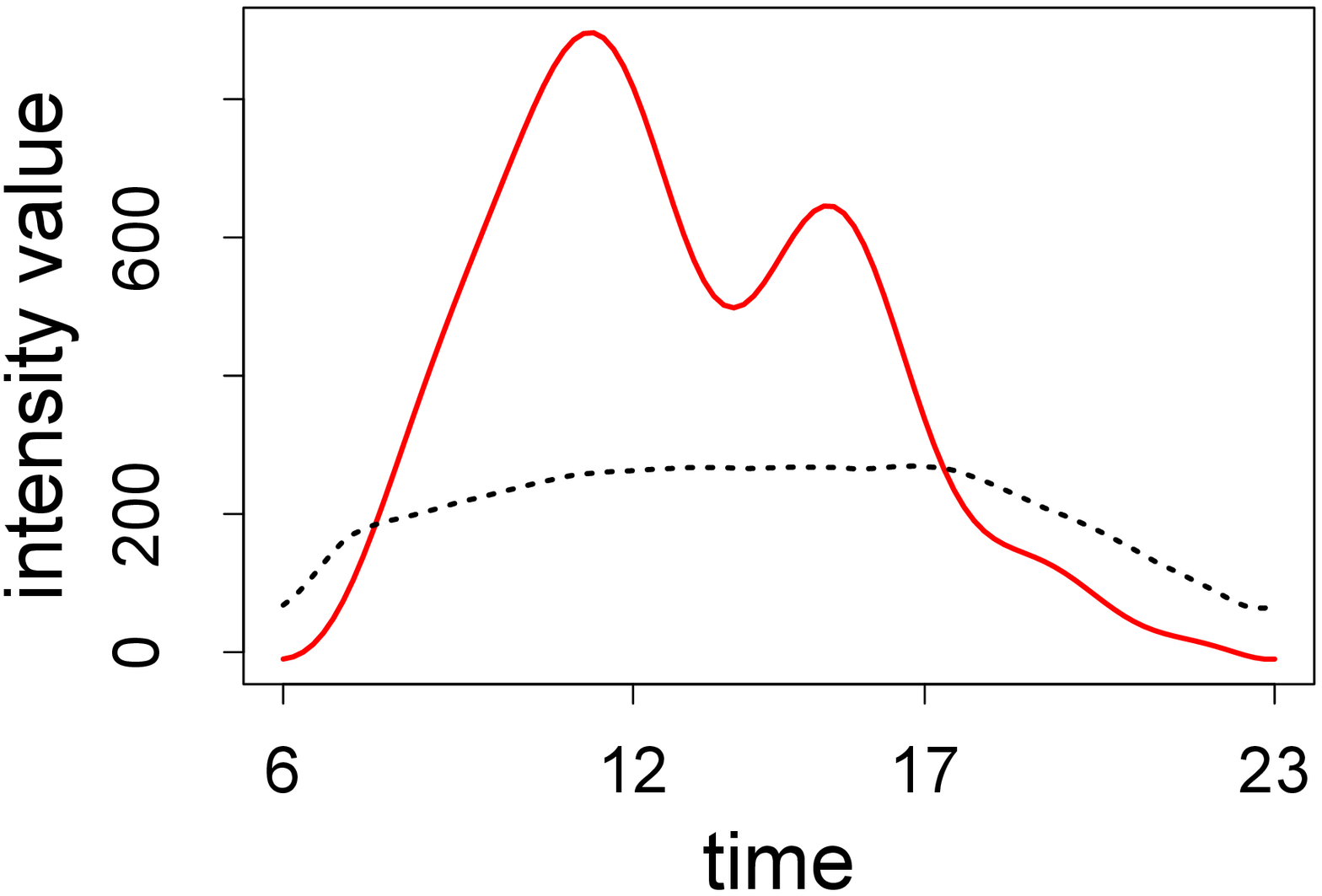}}\hspace{0.2cm}%\\
\subfloat[33166, weekend]{
\includegraphics[scale=0.2]{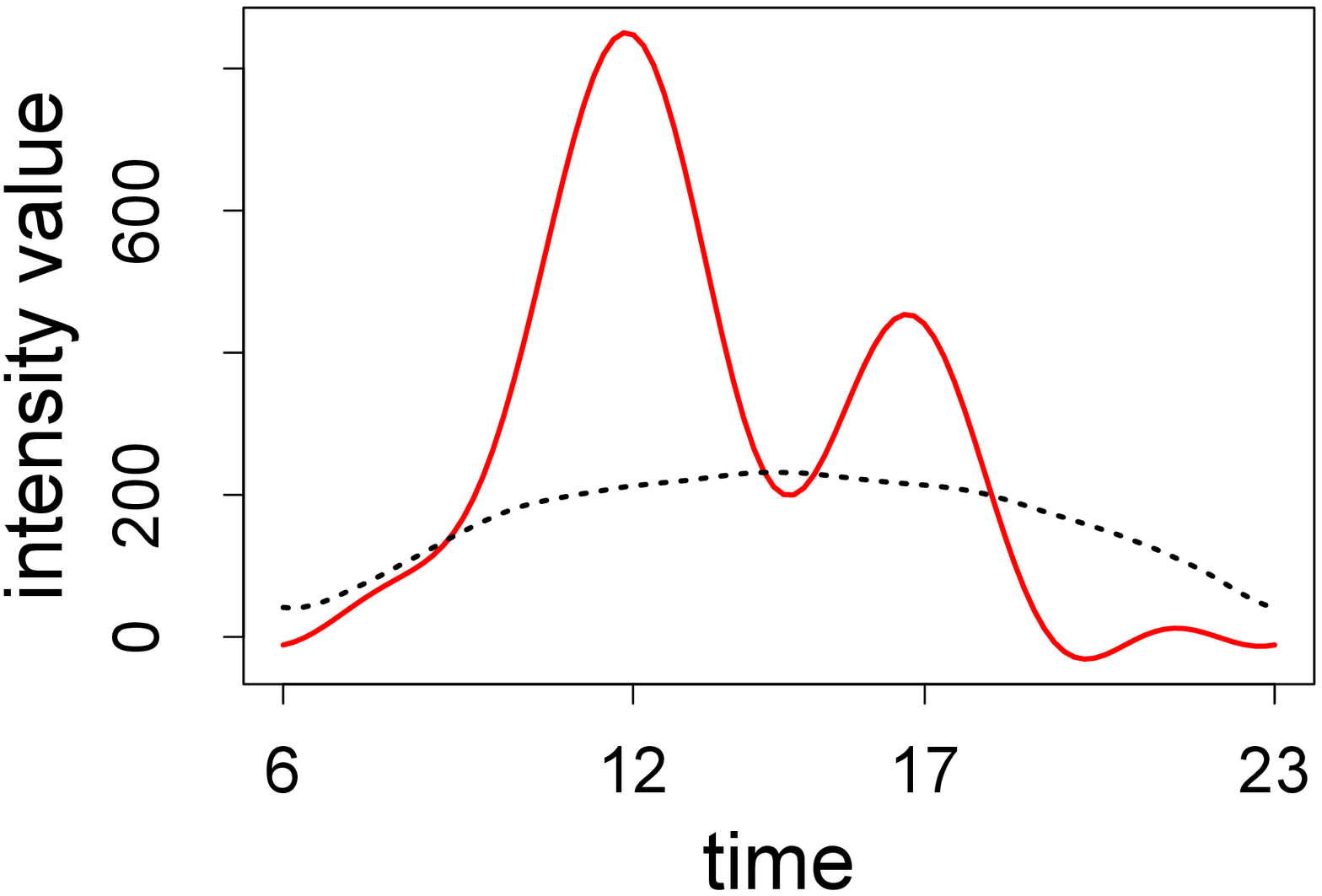}}
\caption{Contrasts between the cohort (black dotted) and individual (red solid)  physical activity intensity profiles for the participants 31212 and 33166.}
\label{fig:individual-profile}
\end{figure}

For each individual, the local stream, when fed into the proposed algorithm, also yields a weekday individual dynamic activity profile and a weekend one. Shown in Figure \ref{fig:individual-profile} are example snapshots of such profiles overlaid by the corresponding cohort profiles. Each panel of Figure \ref{fig:individual-profile} then provides the contrast between the cohort profile and an individual profile. For instance, in a real fitness app, the individual  identified by the SEQN (respondent sequence number) 31212 (age 23.4)  may be presented with Figures \ref{fig:individual-profile}(a)(b) and knows that his/her activity profile is overall much lower than the cohort average. In contrast, the individual 33166 (age 27.5) can learn from Figures \ref{fig:individual-profile}(c)(d) that he/she has an overall higher activity profile relative to the cohort one.

\begin{figure}
\centering
\subfloat[time = 6:00]{
\includegraphics[scale=0.2]{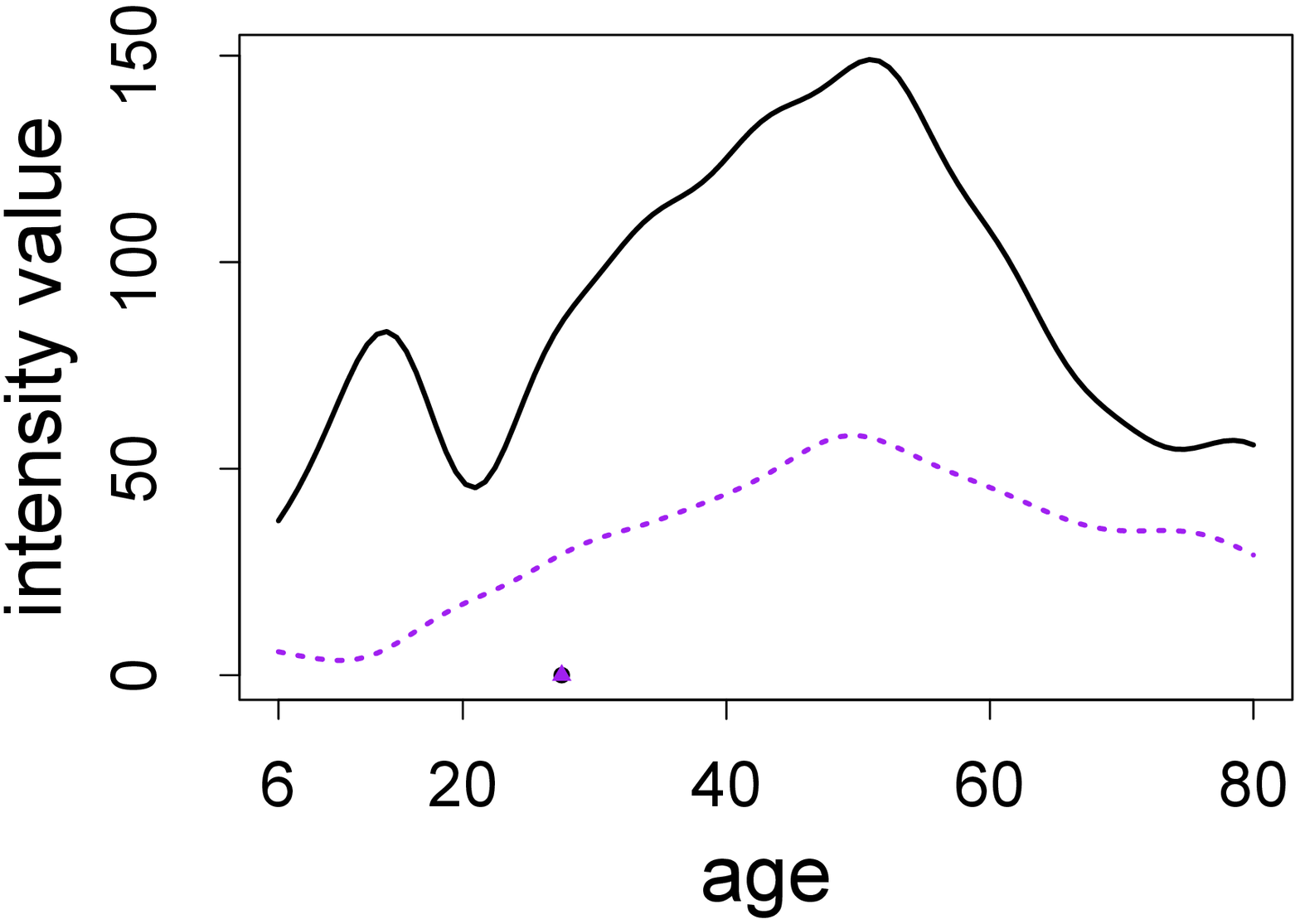}}\hspace{0.3cm}
\subfloat[time = 12:00]{
\includegraphics[scale=0.2]{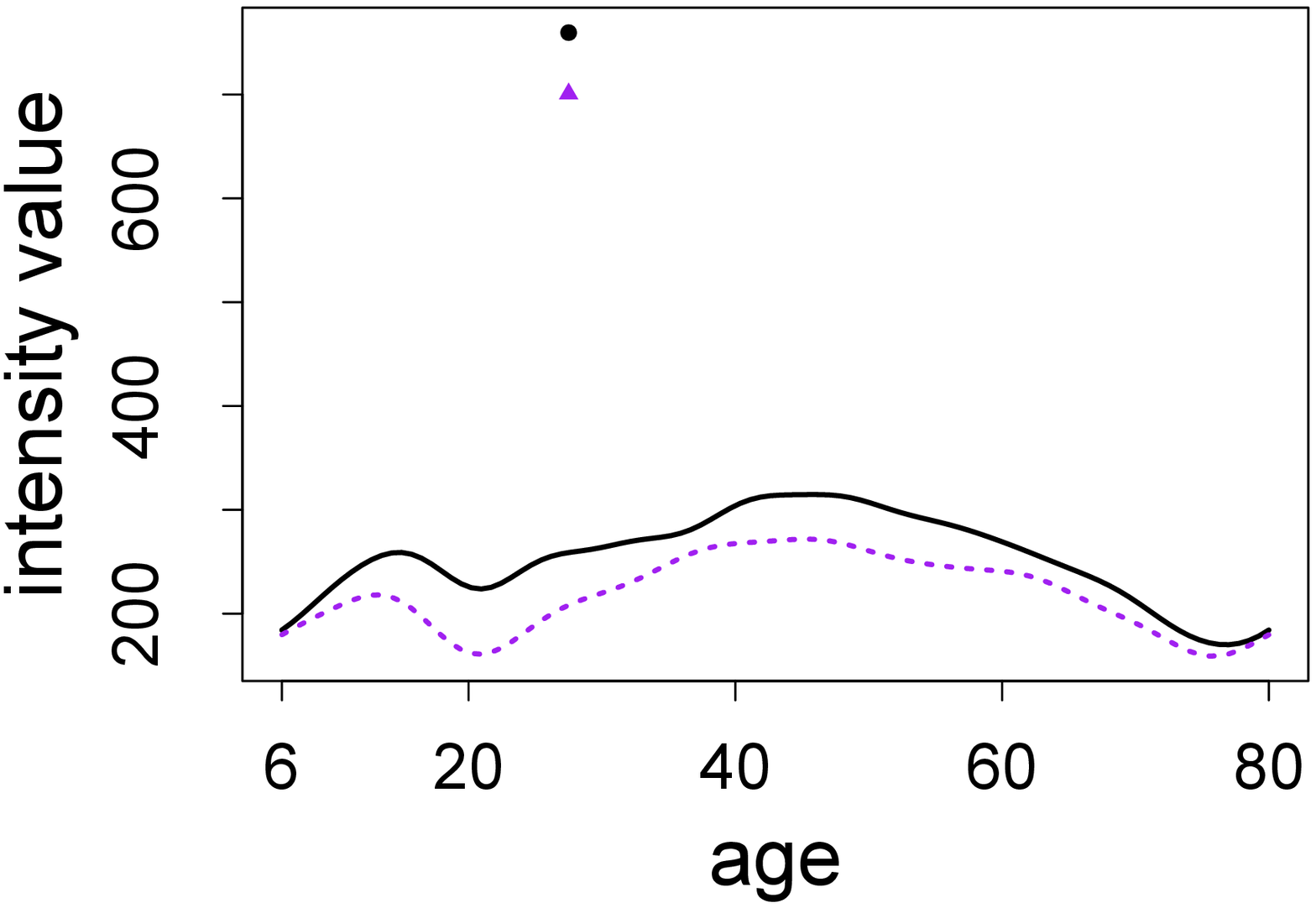}}\hspace{0.3cm}
\subfloat[time = 17:00]{
\includegraphics[scale=0.2]{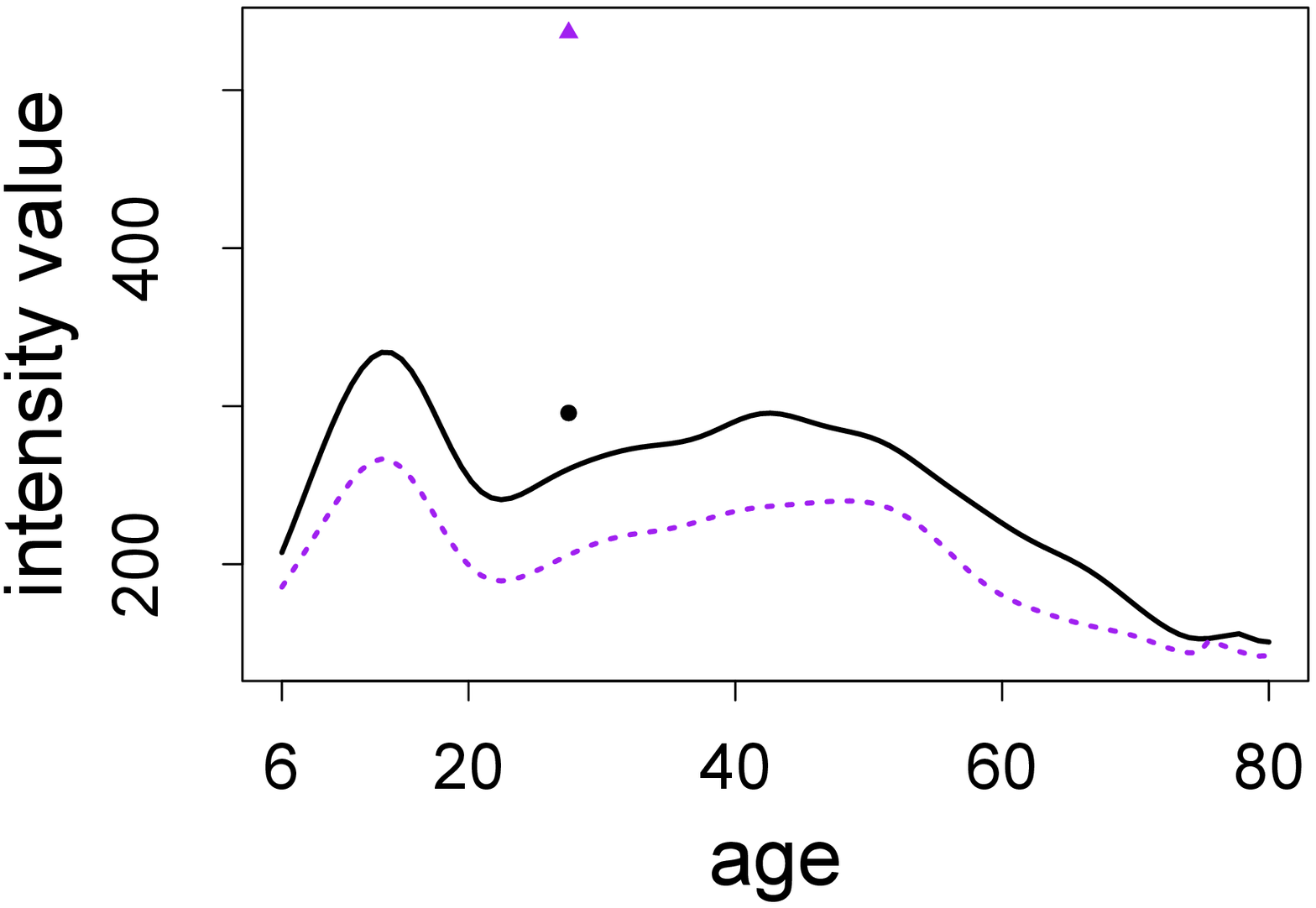}}\hspace{0.3cm}
\subfloat[time = 23:00]{
\includegraphics[scale=0.2]{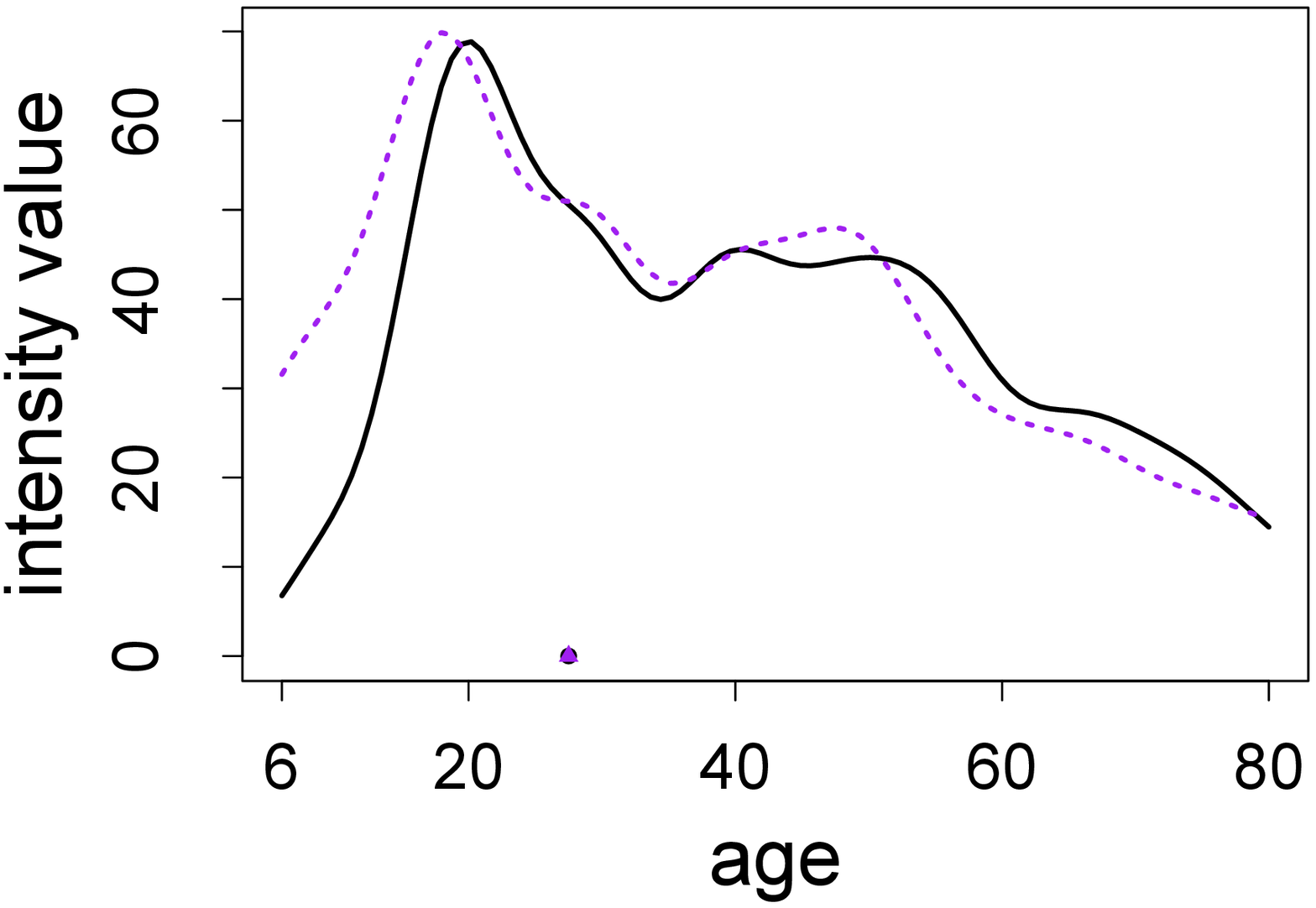}}
\caption{Estimated weekday (black solid) and weekend (purple dotted) mean activity intensity curves of age at different time points of the day, overlaid by the weekday (black circle) and weekend (purple triangle) mean intensity of the individual 33166 at the corresponding time points.}
\label{fig:age-profile}
\end{figure}

It may be also informative for an individual to know how physical activity intensity level changes with age. For this purpose, we  consider an alternative global stream formed by the pairs $(A_i,Y_i)$ collected from all participants, where $A_i$ is the age of an individual and $Y_i$ is the recorded intensity level at a specific time $t\in[6,23]$.  Applying the proposed algorithm to this data stream yields a dynamic activity profile as a function of the age at the given $t$. Examples of such profile are shown in Figure \ref{fig:age-profile}, overlaid by the average intensity level of the individual 33166. In a real fitness app, this individual can then learn from these plots that his/her activity intensity is overall significantly higher than the average activity intensity of the cohort of {age 27.5} at the time $t$=12:00 and $t$=17:00, especially during the weekend. In addition, the nearly zero activity levels of this individual at the time $t$=23:00 and $t$=6:00 reveal that the individual is likely sleeping during 23:00--6:00, indicating a regular circadian rhythm of the individual.  Figure \ref{fig:age-profile}(a) also suggests that children tend to have more time in bed in the weekend morning, while most middle-aged adults get up at or before 6:00 maybe for housework. Moreover, Figure \ref{fig:age-profile}(d) shows that late adolescents and young adults tend to enjoy more late-night activities such as parties around 23:00 especially during the weekend.

\section{Discussion}\label{sec:discussion}
Nonparametric regression for streaming data has been investigated under the framework of regret theory \citep[e.g.,][and references therein]{rakhlin2014} and relative efficiency \citep{yang2021online}, where the regret framework evaluates the performance of an estimator/algorithm against a benchmark class of functions, while the relative efficiency concerns how well a streaming algorithm performs relative to its non-streaming counterpart. It is possible to consider space complexity under the regret framework, for instance, by replacing the $\|\cdot\|_{L^2}$ or $\|\cdot\|_\infty$ loss functions in Section \ref{subsec:space-np-reg} with the worst-case regret. 
It is then of interest to investigate tight lower bounds on memory footprints of a class of streaming estimators that share a common upper bound on the max regret.  
For the framework of relative efficiency, a similar concept of space complexity and minimaxity may be introduced to quantify the minimal amount of memory required by any one-pass counterpart of a non-streaming estimator as a function of both the sample size and the relative efficiency.  
Studies on the interplay between statistical efficiency and memory usage under the framework of regret theory or relative efficiency, requiring substantial development of new techniques and theories, are beyond the scope of our paper and thus left for future investigation.

\bigskip
\begin{center}
	{\large\bf SUPPLEMENTARY MATERIAL}
\end{center}

	The supplementary document contains additional numeric experiments and results, an algorithmic description of the proposed one-pass estimator, a fully data-driven approach for selecting tuning parameters, a brief description on Fourier extension, and proofs for the theoretical results in Section \ref{sec:convergence rate}.

\bibliographystyle{apalike}
\bibliography{ref}

%\includepdf[pages=-]{supplement.pdf}

\end{document}